# EFFICIENT IMPORTANCE SAMPLING FOR BINARY CONTINGENCY TABLES[1]

BY JOSE H. BLANCHET

*Columbia University*

Importance sampling has been reported to produce algorithms with excellent empirical performance in counting problems. However, the theoretical support for its efficiency in these applications has been very limited. In this paper, we propose a methodology that can be used to design efficient importance sampling algorithms for counting and test their efficiency rigorously. We apply our techniques after transforming the problem into a rare-event simulation problem—thereby connecting complexity analysis of counting problems with efficiency in the context of rare-event simulation. As an illustration of our approach, we consider the problem of counting the number of binary tables with fixed column and row sums, $c_j$'s and $r_i$'s, respectively, and total marginal sums $d = \sum_j c_j$. Assuming that $\max_j c_j = o(d^{1/2})$, $\sum c_j^2 = O(d)$ and the $r_j$'s are bounded, we show that a suitable importance sampling algorithm, proposed by Chen et al. [*J. Amer. Statist. Assoc.* **100** (2005) 109–120], requires $O(d^3\varepsilon^{-2}\delta^{-1})$ operations to produce an estimate that has $\varepsilon$-relative error with probability $1-\delta$. In addition, if $\max_j c_j = o(d^{1/4-\delta_0})$ for some $\delta_0 > 0$, the same coverage can be guaranteed with $O(d^3\varepsilon^{-2}\log(\delta^{-1}))$ operations.

**1. Introduction.** We are interested in the complexity analysis of sequential or state-dependent importance sampling algorithms (SIS) for counting problems. The development of algorithms for approximate counting in polynomial time has been a topic of great interest in theoretical computer science [see Valiant (1979)]. Successful techniques have been developed for efficient approximate counting based on the Markov Chain Monte Carlo (MCMC)

Received December 2006; revised August 2008.
[1]Supported in part by NSF Grant DMS-05-95595.

method [see the texts by Sinclair (1993) and Jerrum (2003) for detailed information on these techniques]. A different class of randomized algorithms for approximate counting, based on importance sampling, has received substantial attention recently [basic notions on importance sampling are discussed in Section 2.2; for additional background on importance sampling, see Asmussen and Glynn (2007) and Liu (2001)]. Chen et al. (2005) proposed an algorithm based on importance sampling for counting the number of bipartite graphs with a given degree sequence. They tested their algorithm empirically and observed that it achieved excellent performance. Recently, Blitzstein and Diaconis (2008) have also used importance sampling algorithms for approximately counting the number of acyclic and undirected graphs with a given degree sequence. In addition, Rubinstein (2007) and Botev and Kroese (2008) have applied adaptive importance sampling algorithms to a variety of combinatorial problems, including counting and optimization. Although many of these algorithms based on importance sampling seem to have excellent practical performance, the theoretical framework to carry through a rigorous analysis of their performance is still under development.

Our purpose is to illustrate a framework that can be used to design efficient importance sampling algorithms for counting and provide a rigorous analysis of their computational complexity. Our method provides a direct connection between asymptotic approximations and efficient importance sampling, and we believe that the principle underlying this connection can be applied in substantial generality. In order to illustrate our proposed techniques, we shall consider the problem of counting the number of 0–1 matrices with specified column and row sums—these types of matrices are also called binary contingency tables in statistical applications. In the context of graph theory, this problem is equivalent to that of counting the number of bipartite graphs with a given degree sequence.

Returning to the problem that we consider here, we mention that statistical analysis of binary contingency tables is a problem that has been motivated by several application domains, including some in biology as explained in Chen et al. (2005). Our goal is to provide rigorous support for the observed experimental efficiency of a class of SIS algorithms proposed by Chen et al. (2005) for counting binary contingency tables. Formally, the problem consists in developing fast computational algorithms for counting the number of solutions $\{x_{ij} : 1 \leq i \leq m, 1 \leq j \leq n\}$ to

$$\sum_{j=1}^{n} x_{ij} = r_i, \qquad i \in \{1, 2, \ldots, m\}, \tag{1}$$

$$\sum_{i=1}^{m} x_{ij} = c_j, \qquad j \in \{1, 2, \ldots, n\}, \tag{2}$$

$$x_{ij} \in \{0,1\}.$$

Let us define $d = \sum_{i=1}^{m} r_i = \sum_{j=1}^{n} c_j$. Our complexity analysis is performed by sending $d \nearrow \infty$ (and $n, m \nearrow \infty$ as well) in the context of sparse matrices under regularity conditions. In particular, we assume that $\max_{j \leq n} c_j = o(d^{1/2})$, $\sum_{j=1}^{n} c_j^2 = O(d)$ and that the $r_i$'s are bounded. We shall construct and analyze a SIS-based estimator, that is, $\varepsilon$-close (in relative terms) to the total number of solutions to (1)–(2) with coverage probability of at least $1 - \delta$ and that requires $O(d^3 \varepsilon^{-2} \delta^{-1})$ operations as $d \nearrow \infty$ and $\varepsilon, \delta \searrow 0$ for its construction. Moreover, by imposing an additional growth condition of the form $\max_{i \leq n} c_j = o(d^{1/4 - \delta_0})$ for some $\delta_0 > 0$ as $d \nearrow \infty$, we obtain that SIS yields a "fully polynomial randomized approximation scheme" in the sense that $O(d^3 \varepsilon^{-2} \log(\delta^{-1}))$ operations are sufficient to produce an estimator that has $\varepsilon$-relative precision with probability $1 - \delta$.

The proposed strategy for constructing and designing SIS algorithms proceeds as follows. The first step is to transform the counting problem into a so-called rare-event estimation problem; that is, the problem of estimating the probability of a rare event. Such probability is given by the ratio of the number of required solutions ($x_{ij}$'s) satisfying both (1) and (2) divided by the number of solutions satisfying only the set of constraints (2) (i.e., there is no restriction on the row sums), which can be easily computed. The second step is to recognize that such probability can be characterized by a system of linear equations, which is obtained by conditioning on the first increment of a suitably defined $m$-dimensional random walk. The solution to this system of equations provides the means of constructing the optimal importance sampling distribution, which is state-dependent. Such optimal importance sampling distribution corresponds to the so-called Doob's $h$-transform which arises in the context of positive harmonic functions. The next step is to use results developed by McKay (1984) [see also Greenhill, McKay and Wang (2006) and Bekessy, Bekessy and Komlos (1972)] that approximate the number of solutions satisfying constraints (1) and (2) in the context of large and sparse matrices that we have adopted here. We then use these approximations (which we have extended in Theorem 1 to cover our assumptions) to construct an importance sampling distribution that mimics the behavior of the optimal importance sampler (thus, the better the approximation the closer the importance sampler to the optimal one).

It turns out that the importance sampling algorithm suggested by the previous strategy coincides with one of the algorithms studied in Chen et al. (2005). Our results imply that in the context of large and sparse matrices satisfying the assumptions indicated above, the variance of the estimator obtained by the procedure has the best possible asymptotic performance. That is, the coefficient of variation of the estimator (the ratio of the standard deviation to the probability of interest) remains bounded as $d \nearrow \infty$.

In the context of rare-event simulation, an estimator that has a bounded coefficient of variation is said to be *strongly efficient* (a notion that will be reviewed in Section 4). Moreover, we show that if $\max_{i \leq n} c_j = o(d^{1/4-\delta_0})$ then the proposed estimator is *exponentially efficient*, a concept that is introduced in Section 4, and in particular implies strong efficiency. In particular, exponential efficiency allows to conclude that under our assumptions, the proposed SIS estimator is $\varepsilon$-close in relative terms with coverage probability $1 - \delta$ and has complexity $O(d^3 \varepsilon^{-2} \log(\delta^{-1}))$ as indicated above, in contrast to a complexity of order $O(d^3 \varepsilon^{-2} \delta^{-1})$ corresponding to strong efficiency.

A recent algorithm for counting binary tables developed by Bezáková, Bhatnagar and Vigoda (2006), based on MCMC techniques (simulated annealing), has been shown to have complexity of (roughly) $O(d^3(nm)^2 \max_{i,j}(c_i, r_j))$ operations. However, it is important to note that the procedure proposed by Bezáková, Bhatnagar and Vigoda (2006) works in complete generality (i.e., it does not require the sparsity assumptions imposed here). Other algorithms based on MCMC have been devised for counting binary contingency tables with certain regularity conditions on the degree sequence (such as ours). For instance, Kim and Vu (2003) assumed that $\max(r_i, c_j) = o(d^{1/3})$ and proposed an algorithm that allows to generate an almost uniform bipartite graph (with given degree sequence) in time $O(d^3)$. Kannan, Tetali and Vempala (1997) also study the problem of uniform generation of bipartite graphs with given degrees.

In a recent paper [Bayati, Kim and Saberi (2007)] used ideas based on SIS to construct an algorithm for generation of simple graphs with a given degree sequence (a slightly different problem than the one that we study here). Under regularity conditions (similar to those imposed here), they proved that their proposed algorithm has excellent performance for asymptotically uniform generation, basically linear complexity, which makes the algorithm optimal in the sense that no faster complexity rate is possible. Their methods seem completely different from those developed here. In particular, we do not require the explicit use of concentration inequalities but instead the application of bounds related to Lyapunov inequalities for Markov chains. In addition, our methods suggest a natural way to develop efficient SIS in a variety of settings—basically if the optimal importance sampling distribution is described by a Markov chain and there are asymptotic approximations for the quantity of interest.

It is important to emphasize that although our complexity analysis of SIS suggests very good performance [which is validated by the computer experiments performed by Chen et al. (2005)], such performance can only be guaranteed under certain regularity conditions. This has been noted by Bezáková et al. (2007) who constructed a counterexample showing that a SIS related to the one presented here [also proposed by Chen et al. (2005)] can have exponential time complexity if the degree sequence is allowed to grow arbitrarily. Nevertheless, one of the main points that we intend to communicate is that the general method outlined here could be adapted to specific contexts in which the problem at hand seems to have certain regularity properties that allow to develop approximations. The strategy would be then to enhance the approximations by means of efficient computational algorithms that can be shown to have desirable complexity properties in an asymptotic regime related to the developed approximations.

The basic principles behind the design and analysis of the SIS algorithms discussed here can be applied more broadly. For instance, Blanchet and Glynn (2008) apply these principles in the context of first-passage time probabilities with an emphasis on one-dimensional random walk problems with general heavy-tailed increments (which are of particular interest in insurance and queueing). Another example is given in Blanchet and Liu (2008), which develops strongly efficient rare-event simulation algorithms for large deviation probabilities of regularly varying random walks. The analysis of SIS algorithms in rare-event simulation involves constructing so-called Lyapunov functions, which are solutions to certain inequalities that are used in stability analysis of Markov processes. The use of Lyapunov inequalities in the context of counting problems as the one considered here is particularly interesting because the dimension of the state-variables of the underlying Markov process (in our context $m$) is growing. As we shall discuss in Section 5, the construction of a suitable Lyapunov function often requires a good understanding of the local likelihood ratio obtained at each step of the simulation.

The rest of the paper is organized as follows. In Section 2, we relate the problem of counting binary contingency tables to its rare-event simulation counterpart. Basic notions involving importance sampling and the optimal change-of-measure are also discussed in Section 2. Section 3 develops asymptotic approximations that are then used in the construction of the algorithm that we analyze. Section 4 introduces efficiency notions that are applied in rare-event simulation and discusses connections to related ideas used in the context of approximate counting. The complexity analysis of the counting algorithm, which leads to the proof of our main result, namely Theorem 2, is given in Section 5.

**2. Counting, rare-event simulation and importance sampling.** A 0–1 table with specified marginals is a binary array (0–1 elements) of dimensions $m \times n$ such that the sum of the elements in the $i$th row equals $r_i$ ($i \in \{1, \ldots, m\}$) and the corresponding sum over the $j$th column equals $c_j$, $j \in \{1, \ldots, n\}$.

*Notational convention.* Throughout the rest of the paper, we shall use the notation $\mathbf{c} = (c_1, \ldots, c_n)$, $\mathbf{r} = (r_1, \ldots, r_m)$ and $\sum_{j=1}^n c_j = d = \sum_{i=1}^m r_i$. In addition, we shall reserve the use of boldface letters to denote vectors or high-dimensional objects. Random variables are denoted using capital letters and the use of lower case is restricted to deterministic quantities (including specific realizations of random objects). We also use the notation $f(t) = O(g(t))$ if there exists a constant $m_1 \in (0, \infty)$ such that $|f(t)| \leq m_1 g(t)$; if, in addition, $|f(t)| \geq m_2 g(t)$ for some $m_2 \in (0, \infty)$, then $f(t) = \Theta(g(t))$. Finally, we say that $f(t) = o(g(t))$ as $t \nearrow \infty$ if $f(t)/g(t) \longrightarrow 0$ as $t \nearrow \infty$.

2.1. *Counting binary tables via a random walk computation.* We are interested in developing an importance sampling algorithm that allows to efficiently count the number of such arrays, which we shall denote by $\mu(\mathbf{r}, \mathbf{c})$. First, note that the number of tables with $m$ rows and given only column marginals $\mathbf{c}$ is

$$\eta(\mathbf{c}, m) \triangleq \binom{m}{c_1} \cdots \binom{m}{c_n}.$$

So, the number of tables with given column and row marginals, $\mathbf{c}$ and $\mathbf{r}$, respectively, can be evaluated via

$$\eta(\mathbf{c}, m) \cdot P(\mathbf{T}(\mathbf{Y}) = \mathbf{r}),$$

where $\mathbf{T}(\mathbf{Y}) \in \mathbb{R}^m$ is the row marginals of a table $\mathbf{Y}$ sampled uniformly over the space of tables with column marginals $\mathbf{c}$ [i.e., $\mathbf{T}(\mathbf{Y})_i = \sum_{j=1}^n \mathbf{Y}_{i,j}$ for $1 \leq i \leq m$]. As a consequence, the problem of counting the number of binary contingency tables is equivalent to that of estimating $P(\mathbf{T}(\mathbf{Y}) = \mathbf{r})$. We can see that efficient estimation of this probability with good relative precision is not straightforward because the probability in question may become arbitrarily small as the size of the table increases. In other words, the event $\{\mathbf{T}(\mathbf{Y}) = \mathbf{r}\}$ would typically be rare.

We shall formulate the problem of estimating $P(\mathbf{T}(\mathbf{Y}) = \mathbf{r})$ as a sequential rare-event simulation problem involving a suitably defined random walk. Define an $m$-dimensional random walk (rw) $\mathbf{S} = (\mathbf{S}_k : 0 \leq k \leq n)$ as follows. Given vectors of nonnegative integers $\mathbf{c}$ and $\mathbf{r}$, set $\mathbf{S}_0 = \mathbf{r}$ and (for $k \in \{1, \ldots, n\}$) define $\mathbf{S}_k = \mathbf{S}_{k-1} - \mathbf{X}_k$ where $\mathbf{X}_k$ is a 0–1 entry vector (of dimension $m$) with uniform distribution over the space of configurations $(x_{k,1}, \ldots, x_{k,m})$ such that $\sum_{j=1}^m x_{k,j} = c_k$ and (for $1 \leq j \leq m$) $x_{k,j} \in \{0, 1\}$. The vector $\mathbf{X}_k$ represent the $k$th column of the table. The random vectors $(\mathbf{X}_k : 1 \leq k \leq n)$ are assumed to be independent. Finally, let us write $P_{\mathbf{r}, \mathbf{c}}(\cdot)$ for the probability law generated by the rw $\mathbf{S}$ subject to $\mathbf{S}_0 = \mathbf{r}$ and $E_{\mathbf{r}, \mathbf{c}}(\cdot)$ for the corresponding expectation operator. Note that $P_{\mathbf{r}, \mathbf{c}}(\cdot)$ is defined via a time inhomogeneous Markov chain; this is because the distributions of the $\mathbf{X}_k$'s change in time according to $\mathbf{c}$.

Observe that

$$u(\mathbf{r}, \mathbf{c}) \triangleq P_{\mathbf{r}, \mathbf{c}}(\mathbf{S}_n = 0) = P(\mathbf{T}(\mathbf{Y}) = \mathbf{r}).$$

As we shall see, the function $u(\cdot)$ can be used to describe the conditional distribution of the rw $\mathbf{S}$ given that $\mathbf{S}_n = 0$. In turn, such description is highly relevant for the design of efficient importance sampling algorithms.

2.2. *Basic notions on importance sampling.* We shall briefly discuss basic concepts related to importance sampling and then apply these concepts to the random walk computation described in the previous section. For a more detailed discussion on importance sampling, see, for instance, Asmussen and Glynn (2007), Glynn and Iglehart (1989) and Liu (2001).

Suppose that we want to estimate $P(Z \in A) > 0$, for a given random object taking values on a space $\mathcal{X}$ with a $\sigma$-field $\mathcal{B}$. Let us define the probability measure $F_Z(dz)$ on $(\mathcal{X}, \mathcal{B})$ via $F_Z(dz) = P(Z \in dz)$, so that

$$P(Z \in A) = \int_A F_Z(dz).$$

We say that a probability measure $G(dz)$ on $(\mathcal{X}, \mathcal{B})$ is an *admissible* choice for an *importance sampler* or *change-of-measure* if $F_Z(dz)I_A(z)$ is absolutely continuous with respect to $G(dz)$ [note that $F_Z(dz)I_A(z)$ is not necessarily a probability measure]. In other words, if $G(dz)$ is admissible then $G(B) = 0$ implies $P(A \cap B) = 0$, so that the Radon–Nikodym derivative $L(z) = I_A(z)(dF_Z/dG)(z)$ is well defined. If $G(dz)$ is admissible, then

$$(3) \qquad P(Z \in A) = E^G L(Z) = \int L(z) G(dz).$$

Here, we are using $E^G(\cdot)$ to denote an expectation that is computed under the probability measure $G(\cdot)$ [similarly, we will use $\mathrm{Var}_G(\cdot)$ for variances under $G(\cdot)$]. The random variable $L$, which is clearly an unbiased estimator of $P(Z \in A)$, is an *importance sampling estimator*. In some discussions on importance sampling, the likelihood $W(z) \triangleq (dF_Z/dG)(z)$, when is well defined, is said to be the "importance sampling weight" [see, e.g., Liu (2001)]. When $W(z)$ is well defined, then one can write $L(z) = W(z)I_A(z)$.

The idea behind importance sampling is to take advantage of representation (3) in order to estimate $P(Z \in A)$. In particular, one can simulate $k$ i.i.d. (independent and identically distributed) copies of $Z$ using the distribution $G(\cdot)$ and output the estimator

$$(4) \qquad W_{IS}^{(k)} = \frac{1}{k} \sum_{j=1}^{k} L(Z_j).$$

By the LLNs and identity (3), $W_{IS}^{(k)}$ is a consistent estimator of $P(Z \in A)$ as $k \nearrow \infty$. Note that importance sampling can in principle achieve zero variance. Indeed, if one chooses as change-of-measure

$$G^*(dz) = P(Z \in dz | Z \in A) = \frac{P(Z \in dz) I(z \in A)}{P(Z \in A)},$$

we obtain, say if $m = 1$ and $Z = z_1$,

$$W_{IS} = L(z_1) = P(Z \in dz_1) \left[ \frac{P(Z \in dz_1)}{P(Z \in A)} \right]^{-1}$$
$$= P(Z \in A).$$

So, our estimate of $P(Z \in A)$ is exact (in particular, it has zero variance). Obviously, such importance sampling estimator is not feasible to implement in practical cases because it requires knowledge of $P(Z \in A)$, which is the quantity of interest. However, the form of the zero-variance importance sampler indicates that a good change-of-measure should be similar to the conditional distribution of $Z$ given that $Z \in A$.

2.3. *Zero-variance importance sampler for binary contingency tables.* Note, conditioning on $\mathbf{X}_1$, that $u(\mathbf{r}, \mathbf{c})$ satisfies

$$u(\mathbf{r}, \mathbf{c}) = E_{\mathbf{r}, \mathbf{c}}(u(\mathbf{S}_1, \boldsymbol{\rho}_1)),$$

where $\boldsymbol{\rho}_1 = (c_2, \ldots, c_n)$; note that the dimension of $\boldsymbol{\rho}_1$, which we shall denote by $size(\boldsymbol{\rho}_1)$, equals $n - 1$. More generally, at time $0 \le k \le n - 1$

$$(5) \qquad u(\mathbf{s}_k, \boldsymbol{\rho}_k) = E_{\mathbf{s}_k, \boldsymbol{\rho}_k}(u(\mathbf{S}_{k+1}, \boldsymbol{\rho}_{k+1})),$$

where $\boldsymbol{\rho}_{k+1} = (c_{k+2}, \ldots, c_n)$ and $size(\boldsymbol{\rho}_k) = n - k$. If we denote the empty vector by the symbol $*$, we must have that $u(\mathbf{0}, *) = 1$ and $u(\mathbf{r}, *) = 0$ for $\mathbf{r} \ne \mathbf{0}$.

Let us define a Markov kernel $Q^*_{\boldsymbol{\rho}_k}(\cdot)$ (for $0 \le k \le n$) via

$$Q^*_{\boldsymbol{\rho}_k}(\mathbf{s}_k, \mathbf{s}_{k+1}) = \binom{m}{c_{k+1}}^{-1} \frac{u(\mathbf{s}_{k+1}, \boldsymbol{\rho}_{k+1})}{u(\mathbf{s}_k, \boldsymbol{\rho}_k)}.$$

Note that (5) guarantees that $Q^*_{\boldsymbol{\rho}_k}$ is a well defined Markov kernel [i.e., the probabilities $Q^*_{\boldsymbol{\rho}_k}(\mathbf{s}_k, \cdot)$ is a probability mass function]. In order to simplify the notation in what follows, we shall drop the explicit dependence on the subindex $\boldsymbol{\rho}_k$. If one could use $P^{Q^*}(\cdot)$ as importance sampler for simulation [i.e., simulate the process $(\mathbf{S}_k : 0 \le k \le n)$ according to transitions generated by $Q^*(\cdot)$], then our likelihood ratio estimator would be (given $\mathbf{S}_0 = \mathbf{r}$ and $\boldsymbol{\rho}_0 = \mathbf{c}$)

$$I(\mathbf{S}_n = \mathbf{0}) \prod_{k=0}^{n-1} \frac{u(\mathbf{S}_k, \boldsymbol{\rho}_k)}{u(\mathbf{S}_{k+1}, \boldsymbol{\rho}_{k+1})} = \frac{u(\mathbf{S}_0, \boldsymbol{\rho}_0)}{u(\mathbf{0}, *)} = u(\mathbf{r}, \mathbf{c}),$$

which has zero variance. Therefore, $Q^*(\cdot)$ corresponds to the zero-variance importance sampling distribution. The kernel $Q^*(\cdot)$ is the so-called Doob's $h$-transform and describes the conditional distribution of the process $\mathbf{S}$ given that $\mathbf{S}_n = \mathbf{0}$; see Doob (1957).

The efficient design of importance sampling algorithms should take advantage of any available information about $u(\cdot)$. For instance, if we know that $u(\cdot)$ is in some sense close to some computable function $v(\cdot)$, then given our previous discussion, it is natural to consider a transition kernel of the form

$$Q(\mathbf{s}_k, \mathbf{s}_{k+1}) = \binom{m}{c_{k+1}}^{-1} \frac{v(\mathbf{s}_{k+1}, \boldsymbol{\rho}_{k+1})}{w(\mathbf{s}_k, \boldsymbol{\rho}_k)},$$

where $w(\mathbf{s}_k, \boldsymbol{\rho}_k)$ is the appropriate normalizing constant that makes $Q(\cdot)$ a well defined Markov transition kernel. Once again, for notational simplicity, we suppress the explicit dependence of $\boldsymbol{\rho}_k$ in $Q(\cdot)$ but keep in mind that $Q(\cdot)$ is a time inhomogeneous Markov kernel. This is the strategy that we shall pursue in the next section in order to describe an importance sampling scheme that can be rigorously shown to be efficient in a context of sparse tables. An early reference that explores the connection between $h$-transforms and importance sampling is Glynn and Iglehart (1989) [see also Asmussen and Glynn (2007) and Juneja and Shahabuddin (2006) for more applications of this idea].

**3. Approximating the optimal change-of-measure and algorithm design.**
In order to apply the strategy outlined at the end of the previous section, we need to find a suitable approximation $v(\cdot)$ to $u(\cdot)$. Results from McKay (1984) and Greenhill, McKay and Wang (2006) will allow us to obtain valuable information on $u(\cdot)$ that we will exploit in order to design an efficient importance sampling algorithm. In order to develop the required approximations, it is useful to introduce some notation.

As we indicated at the beginning of the previous section, $\mu(\mathbf{r}, \mathbf{c})$ represents the number of tables with fixed column sums vector $\mathbf{c}$ and marginal row sums given by $\mathbf{r}$. Note that with this definition of $\mu(\mathbf{r}, \mathbf{c})$ we have $u(\mathbf{r}, \mathbf{c}) = \mu(\mathbf{r}, \mathbf{c})/\eta(\mathbf{c}, m)$. We shall assume that the $c_j$'s are ordered in a nonincreasing way so that $c_1 \geq c_2 \geq \cdots \geq c_n$. Having the $c_j$'s ordered in this way does not affect the asymptotic approximations that we are about to describe, but as we shall see, the ordering is important for the good performance of SIS.

We now introduce some convenient notation as in Greenhill, McKay and Wang (2006) that will be useful throughout the rest of the paper. Given a number $s$ and an integer $k \geq 0$, we define $[s]_k = s(s-1)\cdots(s-k+1)$ and $[s]_0 = 1$. Given a vector $\mathbf{s} = (s_1, \ldots, s_{n_0})$ of dimension $n_0$, we set $[\mathbf{s}]_0 = 1$ and define, for any integer $k \geq 1$,

$$[\mathbf{s}]_k = \sum_{j=1}^{n_0} [s_j]_k, \qquad [\mathbf{s}^k]_1 = \sum_{j=1}^{n_0} s_j^k$$

and also $\mathbf{s}! = s_1! s_2! \cdots s_{n_0}!$.

Define

$$\varphi(\mathbf{r}, \mathbf{c}) = \frac{[\mathbf{c}]_1!}{\mathbf{r}!\mathbf{c}!} \quad \text{and} \quad \alpha(\mathbf{r}, \mathbf{c}) = \frac{[\mathbf{c}]_2[\mathbf{r}]_2}{2[\mathbf{c}]_1^2}.$$

We now are ready state the following result, the proof of which is given at the end of the section. The next theorem is basically an adaptation of results from McKay (1984) [see also Theorem 1.1 of Greenhill, McKay and Wang (2006)].

THEOREM 1. *Assume that* $\max([\mathbf{c}^2]_1, [\mathbf{r}^2]_1) = O(d)$ *and* $\max_{j \leq n, i \leq m}(c_j, r_i) = o(d^{1/2})$ *as* $d \nearrow \infty$. *Then*

$$\mu(\mathbf{r}, \mathbf{c}) \sim \varphi(\mathbf{r}, \mathbf{c}) \exp(-\alpha(\mathbf{r}, \mathbf{c}))$$

*as* $d \nearrow \infty$.

The previous result is slightly different from that of McKay (1984) who required $\max_{i \leq m, j \leq n}\{r_i, c_j\} = o(d^{1/4})$ but did not assume $\max([\mathbf{c}^2]_1, [\mathbf{r}^2]_1) = O(d)$. Further refinements have been given in Theorem 1.3 of Greenhill, McKay and Wang (2006) who introduce additional correction terms by assuming $\max_{i \leq m} r_i \times \max_{j \leq n} c_j = O(d^{2/3})$.

Continuing in the spirit of our discussion at the end of the previous section. We are interested in proposing a function $v(\cdot)$ that mimics the behavior of $u(\cdot)$ in some sense in order to construct our importance sampling algorithm. Theorem 1 suggest using the approximation

$$v(\mathbf{r}, \mathbf{c}) \triangleq \frac{\varphi(\mathbf{r}, \mathbf{c}) \exp(-\alpha(\mathbf{r}, \mathbf{c}))}{\eta(\mathbf{c}, m)}.$$

Let us define $v(\mathbf{0}, *) = 1$ and $v(\mathbf{s}, \boldsymbol{\rho}) = 0$ if at least one component of $\mathbf{s}$ is negative. As we indicated at the end of last section, our discussion of the zero-variance change-of-measure, $Q^*$, suggests designing the importance sampling distribution via a Markov transition kernel of the form

$$(6) \qquad Q(\mathbf{s}_k, \mathbf{s}_{k+1}) = \binom{m}{c_{k+1}}^{-1} \frac{v(\mathbf{s}_{k+1}, \boldsymbol{\rho}_{k+1})}{w(\mathbf{s}_k, \boldsymbol{\rho}_k)},$$

where $\boldsymbol{\rho}_k = (c_{k+1}, \ldots, c_n)$ and

$$w(\mathbf{s}_k, \boldsymbol{\rho}_k) = \sum_{(\mathbf{s}_k, \boldsymbol{\rho}_k) \to (\mathbf{s}_{k+1}, \boldsymbol{\rho}_{k+1})} \binom{m}{c_{k+1}}^{-1} \frac{\varphi(\mathbf{s}_{k+1}, \boldsymbol{\rho}_{k+1})}{\eta(\boldsymbol{\rho}_{k+1}, m)}$$

is the normalizing constant that makes $Q(\cdot)$ a well defined Markov transition kernel. In the previous display and in the discussion that follows, we use $(\mathbf{s}_k, \boldsymbol{\rho}_k) \to (\mathbf{s}_{k+1}, \boldsymbol{\rho}_{k+1})$ to denote an admissible transition step [i.e., $\mathbf{s}_k - \mathbf{s}_{k+1}$ is an $m$-dimensional 0–1 whose components add up to $c_{k+1}$ and $\boldsymbol{\rho}_{k+1} = (c_{k+2}, \ldots, c_n)$].

We shall mention how to simulate transitions under $Q(\cdot)$ right after the precise description of the proposed algorithm below. We will use $P_{\mathbf{r},\mathbf{c}}^Q(\cdot)$ to denote the probability measure induced by the random walk $\mathbf{S}$ under the transition kernel given that $\mathbf{S}_0 = \mathbf{r}$ and $E_{\mathbf{r},\mathbf{c}}^Q(\cdot)$ to denote the corresponding expectation operator associated to $P_{\mathbf{r},\mathbf{c}}^Q(\cdot)$.

Note that under the change-of-measure $P_{\mathbf{r},\mathbf{c}}^Q(\cdot)$ we may have $P_{\mathbf{r},\mathbf{c}}^Q(\mathbf{S}_n \neq 0) > 0$. Therefore, when running an importance sampling algorithm based on transitions according to $Q(\cdot)$ we may obtain realizations for which $\{\mathbf{S}_n \neq 0\}$. A sufficient condition that implies $\{\mathbf{S}_n \neq 0\}$ and which can be easily checked at a time $k < n$ is that the number of strictly positive components of $\mathbf{S}_k$ is less than $c_{k+1}$. So, the path generation under $P_{\mathbf{r},\mathbf{c}}^Q(\cdot)$ will be done sequentially according to the transition kernel $Q(\cdot)$ until time $n$ (in which case we have that the event $\{\mathbf{S}_n = 0\}$ has occurred) or up to the first time $k$ such that the number of strictly positive components of $\mathbf{S}_k$ is less than $c_{k+1}$ (in which case we have that $\{\mathbf{S}_n \neq 0\}$).

In order to explain this path generation scheme more formally, let us define

$$\Phi(\mathbf{S}_k) = \mathrm{card}\{j : \mathbf{S}_{k,j} > 0\}.$$

That is, $\Phi(\mathbf{S}_k)$ is the number of strictly positive components of $\mathbf{S}_k$. Put $c_{n+1} \triangleq 1$ and define a stopping time $\tau$ via

$$\tau = \inf\{0 \leq k \leq n : \Phi(\mathbf{S}_k) < c_{k+1}\}.$$

Observe that when $\{\tau < n\}$ occurs one of the components of the vector $\mathbf{S}_n$ must be negative and, therefore, $\{\mathbf{S}_n \neq 0\}$. On the other hand, if $\mathbf{S}_\tau = 0$, we must have that $\tau = n$ because the $c_i$'s are strictly positive and $\sum_{j=1}^n c_j = d$. Therefore, we have that $\{\mathbf{S}_n = 0\} = \{\mathbf{S}_\tau = 0\}$, and consequently

$$u(\mathbf{r}, \mathbf{c}) = P_{\mathbf{r},\mathbf{c}}(\mathbf{S}_\tau = 0).$$

The path generation scheme that we described before under the measure $P_{\mathbf{r},\mathbf{c}}^Q(\cdot)$ will be done sequentially up to the stopping time $\tau$. Note that the $k$th column, namely $\mathbf{X}_k$, is generated under $P_{\mathbf{r},\mathbf{c}}^Q(\cdot)$ during the course of the path generation only if $\tau > k - 1$. In turn, $\mathbf{X}_k$ is a binary vector such that the sum of its components equals $c_k$ and $P_{\mathbf{r},\mathbf{c}}^Q(\cdot)$ avoids assigning negative components to the random walk $\mathbf{S}$ and, therefore, generation of increments under $P_{\mathbf{r},\mathbf{c}}^Q(\cdot)$ can be performed up to time $\tau$. If $\tau < n$, then the $\tau$th assignment under $P_{\mathbf{r},\mathbf{c}}^Q(\cdot)$ cannot be done and the estimator is just zero. If $\tau = n$, then the table is constructed satisfying the row and column sums. We then conclude that $P_{\mathbf{r},\mathbf{c}}^Q(\cdot)$ is *admissible* in the sense that it does not assign zero mass to outcomes that are possible under $P_{\mathbf{r},\mathbf{c}}(\cdot)$ and for which $\mathbf{S}_n = 0$. In fact, it turns out that the sequential importance sampling algorithm generated by $Q(\cdot)$ coincides with one of the procedures studied by Chen et al. (2005). In order to see this, note that

$$Q(\mathbf{s}_k, \mathbf{s}_{k+1}) = \binom{m}{c_{k+1}}^{-1} \frac{v(\mathbf{s}_k, \boldsymbol{\rho}_k)}{w(\mathbf{s}_k, \boldsymbol{\rho}_k)} \frac{v(\mathbf{s}_{k+1}, \boldsymbol{\rho}_{k+1})}{v(\mathbf{s}_k, \boldsymbol{\rho}_k)}$$

$$\propto \frac{v(\mathbf{s}_{k+1}, \boldsymbol{\rho}_{k+1})}{v(\mathbf{s}_k, \boldsymbol{\rho}_k)} = \frac{\varphi(\mathbf{s}_{k+1}, \boldsymbol{\rho}_{k+1})\eta(\boldsymbol{\rho}_k, m)}{\varphi(\mathbf{s}_k, \boldsymbol{\rho}_k)\eta(\boldsymbol{\rho}_{k+1}, m)}$$

$$\propto \frac{\varphi(\mathbf{s}_{k+1}, \boldsymbol{\rho}_{k+1})}{\varphi(\mathbf{s}_k, \boldsymbol{\rho}_k)} \propto \frac{\mathbf{s}_k!}{\mathbf{s}_{k+1}!} \exp(\alpha(\mathbf{s}_k, \boldsymbol{\rho}_k) - \alpha(\mathbf{s}_{k+1}, \boldsymbol{\rho}_{k+1}))$$

$$\propto \prod_{j \in \{j : \mathbf{s}_{k+1,j} \neq \mathbf{s}_{k,j}\}} (\mathbf{s}_{k,j} \exp(2\gamma_k \mathbf{s}_{k,j})),$$

where $\gamma_k = \sum_{j=2}^{n-k}(\rho_{k,j}^2 - \rho_{k,j})/(2(d - \rho_{k,1}))$ (note that $\text{card}\{j : \mathbf{s}_{k+1,j} \neq \mathbf{s}_{k,j}\} = \rho_{k,1}$). The proportionality relations are introduced to emphasize the dependence of the transition kernel only on the $k+1$th increment, namely $\mathbf{s}_{k+1} - \mathbf{s}_k$. The last line of the previous display coincides with the description given in page 112 of Chen et al. (2005) [the complete details of the computation corresponding to the difference $\alpha(\mathbf{s}_k, \boldsymbol{\rho}_k) - \alpha(\mathbf{s}_{k+1}, \boldsymbol{\rho}_{k+1})$ are given in Section 5; see (17)].

The precise form of the algorithm that we analyze, based on the transition kernel $Q(\cdot)$ defined in (6), is given next.

ALGORITHM 1.

STEP 1. Order the $c_i$'s so that $c_1 \geq \cdots \geq c_n$ and set $\mathbf{s} \longleftarrow \mathbf{r}$, $\boldsymbol{\rho} \longleftarrow \mathbf{c}$, $L \longleftarrow 1$ and $l \longleftarrow 0$.

STEP 2. Let $\mathcal{A} = \{i : s_i > 0\}$ and define $m_\mathcal{A} = \text{card}(\mathcal{A})$. Put $c_1 \longleftarrow \boldsymbol{\rho}_1$, $\boldsymbol{\rho} \longleftarrow (\boldsymbol{\rho}_2, \ldots, \boldsymbol{\rho}_{size(\boldsymbol{\rho})})$ and $l \longleftarrow l + 1$. If $m_\mathcal{A} < c_1$, put $L = 0$ and GO TO Step 3. Otherwise, if $l < n$, then evaluate

$$\gamma \longleftarrow [\boldsymbol{\rho}]_2/2[\boldsymbol{\rho}]_1^2$$

else, if $l = n$ set $\gamma = 0$. Sample $(Y_{i_1}, \ldots, Y_{i_{m_\mathcal{A}}})$ according to the distribution

(7)  $$P(Y_{i_1} = y_{i_1}, \ldots, Y_{i_{m_\mathcal{A}}} = y_{i_{m_\mathcal{A}}}) = \frac{1}{w'} \prod_{j=1}^{m_\mathcal{A}} (s_j \exp(2\gamma s_j))^{y_{i_j}},$$

where $\sum_{j=1}^{m_\mathcal{A}} y_{i_j} = c_1$, $y_{i_j} \in \{0, 1\}$ and

$$w' = \sum_{\{(y_{i_1}, \ldots, y_{i_{m_\mathcal{A}}}) : y_{i_1} + \cdots + y_{i_{m_\mathcal{A}}} = c_1\}} \prod_{j=1}^{m_\mathcal{A}} (s_j \exp(2\gamma s_j))^{y_{i_j}}.$$

Then update

$$L \longleftarrow \frac{w'}{\prod_{j=1}^{m_{\mathcal{A}}}(s_j \exp(2\gamma s_j))^{Y_{i_j}}} L$$

and for $j \in \mathcal{A}$ put $\mathbf{s}_j \longleftarrow \mathbf{s}_j - Y_j$.

STEP 3. If $L = 0$ or $n = 0$ output $L$ and STOP, otherwise, GO TO Step 2.

The output of Algorithm 1 is given by

$$L = \prod_{k=0}^{\tau-1} \frac{w(\mathbf{S}_k, \boldsymbol{\rho}_k)}{v(\mathbf{S}_{k+1}, \boldsymbol{\rho}_{k+1})} I(\mathbf{S}_\tau = 0).$$

Equation (7) describes the distribution of $\mathbf{X}_{j+1}$ given $\mathbf{s}_j$ under $Q(\mathbf{s}_j, \cdot)$. Sampling according to such distribution can be done, adopting the terminology used by Chen et al. (2005), using the so-called "drafting method" for Conditional-Poisson (CP) distributions described by (7).

For completeness, we shall review the basic properties of the drafting procedure; our discussion follows [Chen et al. (2005) and Chen, Dempster and Liu (1994)]. Given a distribution of the form

$$(8) \qquad P(Z_1 = z_1, \ldots, Z_m = z_m) = \frac{1}{\widetilde{w}} \prod_{j=1}^{m} w_j^{z_j} I\left(\sum_{j=1}^{m} z_j = c\right),$$

where $w_j > 0$ and $z_j \in \{0, 1\}$ for $1 \leq j \leq m$, the drafting method allows to both, sampling the $Z_i$'s and efficiently computing the normalizing constant $\widetilde{w}$. The drafting method is a sequential procedure that allows to sample $c$ units without replacement from the set $A_m = \{1, 2, \ldots, m\}$; the $i$th unit has a probability proportional to $w_i$. Let $A_k$, $0 \leq k \leq c$ be the set of selected units after $k$ draws, so that $A_0 = \oslash$ and $A_c$ is the final sample to be obtained. At the $k$th step (with $1 \leq k \leq c$), a unit $j \in A_{k-1}^c$ is selected into the sample with probability

$$p(j, A_{k-1}^c) = \frac{\widetilde{w}(c-k, A_{k-1}^c \setminus \{j\}) w_j}{(c-k+1)\widetilde{w}(c-k+1, A_{k-1}^c)},$$

where

$$\widetilde{w}(i, A) = \sum_{C \subseteq A, \mathrm{card}(C)=i} \left(\prod_{i \in C} w_i\right),$$

$\widetilde{w}(0, A) = 1$ for all $A \subseteq A_m$ and $\widetilde{w}(i, A) = 0$ for $i > \mathrm{card}(A)$. The computation of the $\widetilde{w}(i, A)$'s is performed using the recursion

$$\widetilde{w}(i, A) = \widetilde{w}(i, A \setminus \{j\}) + \widetilde{w}(i-1, A \setminus \{j\}) w_j.$$

For instance, to compute $\widetilde{w}(c, A_m)$, we apply the recursion

$$\widetilde{w}(i, A_j) = \widetilde{w}(i, A_j \setminus \{j\}) + \widetilde{w}(i-1, A_j \setminus \{j\})w_j$$

for $1 \leq i \leq c$ and $i \leq j \leq m$. It follows that computing $\widetilde{w}(c, A_m)$ takes $O(cm)$ operations. Evaluating $p(j, A_m) = p(j, A_0^c)$ then takes $O(cm^2)$ operations. Each of the $p(j, A_k^c)$'s can be evaluated similarly, however, it is more convenient to use Lemma 1 of Chen, Dempster and Liu (1994), which states that

$$p(j, A_k^c) = \frac{w_{i_k} p(j, A_{k-1}^c) - w_j p(j, A_{k-1}^c)}{(c-k)(w_{i_k} - w_j) p(i_k, A_{k-1}^c)}$$

for $1 \leq k \leq c-1$ and $j \in A_k^c$, where $i_k$ is the element selected in the $k$th iteration of the drafting procedure. Therefore, we conclude that it takes $O(cm^2)$ operations to generate a sample from (8) using the drafting method. Chen and Liu (1997) discuss four additional sampling procedures for CP distributions with similar complexity properties. The previous considerations imply that the computational cost per replication of an importance sampling algorithm based on $Q(\cdot)$ is of order $O(m^2 d + n \log(n))$. The contribution of the term $n \log(n)$ corresponds to ordering the $c_j$'s in Step 1 and computing $\gamma$ in Step 2. Note that subsequent updates of $\gamma$ can be done recursively so there is no need to add an extra factor from the fact that the algorithm goes through Step 2 $n$ times. We summarize these observations in the following lemma.

LEMMA 1. *Algorithm 1 requires $O(m^2 d + n \log(n))$ operations to be executed.*

Chen et al. (2005) also proposed a more refined importance sampling procedure which can be explained as follows. Note that we constructed our importance sampling transition kernel, $Q(\cdot)$, via a suitable approximation $v(\mathbf{r}, \mathbf{c})$ of $u(\mathbf{r}, \mathbf{c})$ that is valid as $d \nearrow \infty$. Furthermore, we introduced additional information into $v(\cdot)$ by defining $v(\mathbf{s}, \boldsymbol{\rho}) = 0$ if $\mathbf{s}$ contains at least one negative component. Intuitively, we could have done even better by setting $v(\mathbf{s}, \boldsymbol{\rho}) = 0$ whenever $u(\mathbf{s}, \boldsymbol{\rho}) = 0$, this is the idea behind the refinement proposed by Chen et al. (2005). One immediate difficulty here is the question of how to easily test the pairs $(\mathbf{s}, \boldsymbol{\rho})$ for which $u(\mathbf{s}, \boldsymbol{\rho}) = 0$. This is achieved by making use of a characterization of so-called graphical sequences (i.e., degree sequences that can give rise to a bipartite graph) in terms of certain constraints that can be easily checked during the course of the simulation. Introducing these types of constraints on the support of $Q$ together with asymptotic approximations may help produce efficient importance sampling estimators (in terms of the discussion given in Section 4). However, in our current context, the vanilla version of the importance sampling procedure,

IMPORTANCE SAMPLING FOR COUNTING 15indicated in Algorithm 1, will already be proved to be efficient in a precise mathematical sense to be described in the next section.

PROOF OF THEOREM 1. We follow closely the steps in the proof of Lemma 2.2 and Theorem 1.3 in Greenhill, McKay and Wang (2006) (GMW). First, we introduce their counting model. We consider a set of $d$ labeled points arranged on $m$ cells, say $\varrho_1, \ldots, \varrho_m$. The cell $\varrho_i$ contains $r_i$ elements. Similarly, we consider another set of $d$ labeled points arranged in $n$ cells denoted by $\Xi_1, \ldots, \Xi_n$ and assume that the $j$th cell, $\Xi_j$, contains $c_j$ elements. We then have $2d$ labeled points in total. A partition of the $2d$ elements into $d$ unordered pairs is called a *pairing*. Each pair is denoted by $e = (\boldsymbol{\rho}, \xi)$ where $\boldsymbol{\rho} \in \varrho_i$ for some $1 \leq i \leq m$ and $\xi \in \Xi_j$ for some $1 \leq j \leq n$. We also write $v(\boldsymbol{\rho})$ to denote the cell corresponding to the point $\boldsymbol{\rho}$ and similarly $v(\xi)$ to denote the cell corresponding to the point $\xi$. A *random pairing* is a pairing that is chosen uniformly at random out of the $d!$ possible pairings. Two pairs are called *parallel* if they involve the same cells. An *error* is an unordered set of two parallel pairs.

It follows easily that the probability of obtaining $l \geq 0$ given pairs occurring in a random pairing is $1/[d]_l$. Let $p_d$ be the probability that no errors occur in a random pairing. As noted by GMW, we have

$$\mu(\mathbf{r}, \mathbf{d})\mathbf{r}!\mathbf{c}! = d! p_d,$$

because (up to a permutation in the labels of the elements in each of the cells) each contingency table corresponds to a pairing that has no errors. Therefore, it suffices to estimate $p_d$ which is done, once again following GMW, using inclusion-exclusion and Bonferroni's inequalities. The inclusion–exclusion development is applied as follows. First, given two different pairs $e$ and $e'$ define $B(e; e')$ to be the set of pairs that contain the particular error $\{e, e'\}$. Note that $p_d = 1 - \overline{p}_d$, where $\overline{p}_d$ is the probability that at least one error occurs in a random pairing. In turn, $\overline{p}_d$ is less or equal to the total contribution corresponding to placements with one error, which we denote by $\overline{b}_d^{(1)}$. In particular, $\overline{b}_d^{(1)} = \sum_{\{e,e'\}} P(B(e, e'))$. More generally, let us define $\overline{b}_d^{(k)}$ as the total contribution in the inclusion-exclusion development corresponding to pairings that contain $k$ errors or more. So, for instance, $\overline{b}_d^{(2)} = \sum_{\{\{e,e'\},\{\tilde{e},\tilde{e}'\}\}} P(B(e, e'), B(\tilde{e}, \tilde{e}'))$, the sum runs over sets of two different errors. We then have that for $k \geq 1$

$$(9) \quad \overline{b}_d^{(1)} - \overline{b}_d^{(2)} + \cdots + \overline{b}_d^{(2k-1)} - \overline{b}_d^{(2k)} \leq \overline{p}_d \leq \overline{b}_d^{(1)} - \overline{b}_d^{(2)} + \cdots + \overline{b}_d^{(2k-1)}.$$

Note that $\overline{b}_d^{(k)}$, $k \geq 2$, can be divided in two parts, namely, one that contains errors that do not have a pair in common, which we denote by $\beta_{d,0}^{(k)}$, and another part that contains errors that have pairs in common, which we denote by $\beta_{d,1}^{(k)}$. We define $\beta_{d,0}^{(1)} = \overline{b}_d^{(1)}$ and note that

$$\overline{b}_d^{(1)} = \frac{1}{d(d-1)}\left(\sum_{i=1}^{n} c_i(c_i-1) \sum_{j=1}^{m} r_j(r_j-1)\right)\bigg/2 = \alpha(\mathbf{r},\mathbf{c}) + o(1)$$

as $d \nearrow \infty$. We claim that for each $k \geq 2$

$$\beta_{d,0}^{(k)} = \frac{\alpha(\mathbf{r},\mathbf{c})^k}{k!} + o(1) \tag{10}$$

as $d \nearrow \infty$. To see this, let us first define $\mathcal{N}_k^0$ as the set of ordered $2k$-tuples of pairs $(e_1, e_2, \ldots, e_k, e_1', \ldots, e_k')$ [with $e_j = (\boldsymbol{\rho}_j, \xi_j)$ and $e_j' = (\boldsymbol{\rho}_j', \xi_j')$] satisfying for each $l \in \{1, \ldots, k\}$, $i_l = i_l'$ and $j_l = j_l'$ where

$$v(\boldsymbol{\rho}_l) = i_l, \qquad v(\boldsymbol{\rho}_l') = i_l',$$
$$v(\xi_l) = j_l, \qquad v(\xi_l') = j_l'$$

with $i_l \neq i_s$ and $j_l \neq j_s$ if $l \neq s$.

Note that

$$\beta_{d,0}^{(k)} = \frac{|\mathcal{N}_k|}{2^k k!} \frac{1}{[d]_{2k}}.$$

We claim that

$$|\mathcal{N}_{k+1}| = |\mathcal{N}_k|\left(\left(\sum_{i=1}^{n} [c_i]_2 \sum_{j=1}^{m} [r_j]_2\right) + o(d^2)\right). \tag{11}$$

To verify this claim let us pick an arbitrary element $(e_1, e_2, \ldots, e_k, e_1', \ldots, e_k') \in \mathcal{N}_k$. We obtain an element of $\mathcal{N}_{k+1}$ by adding two parallel pairs $(e_{k+1}, e_{k+1}')$ so that we obtain $k+1$ errors that do not have pairs in common. This is achieved in

$$\left(\sum_{i=1}^{n} [c_i]_2 - \sum_{l=1}^{k} [c_{i_l}]_2\right)\left(\sum_{j=1}^{m} [r_j]_2 - \sum_{l=1}^{k} [r_{j_l}]_2\right)$$

many ways. Now, since $\max_{j \leq n} c_j = o(d^{1/2})$ we have that

$$\frac{\sum_{l=1}^{k} [c_{i_l}]_2}{d} \leq O\left(\frac{\max_{j \leq n} c_j^2}{d}\right) \longrightarrow 0$$

as $d \nearrow \infty$ (a completely analogous estimate also applied to the sum involving the $r_{j_l}$'s). This implies (11) and as a consequence (10). To study $\beta_{d,1}^{(k)}$, it suffices to perform a very rough analysis. Indeed, note that

$$\beta_{d,1}^{(k)} = \sum_{l=2}^{k} O\left(\frac{\sum_{i=1}^{n} c_i^{l+1} \sum_{j=1}^{m} r_j^{l+1}}{d^{l+1}}\right),$$

where the $l$th term in the previous sum corresponds to collections of $k$ errors in which all the pairs belong to $l$ errors or less. Now, we have that

$$\frac{\sum_{i=1}^{n} c_i^{l+1} \sum_{j=1}^{m} r_j^{l+1}}{d^{l+1}} \leq \max_{i \leq n} \frac{c_i}{d^{1/2}} \times \max_{j \leq m} \frac{r_j}{d^{1/2}} \sum_{i=1}^{n} \sum_{j=1}^{m} \frac{c_i^l r_j^l}{d^l}.$$

Since $l \geq 2$ we have that $\sum_{i=1}^{n} \sum_{j=1}^{m} c_i^l r_j^l / d^l = O(1)$ as $d \nearrow \infty$ and we conclude that $\beta_{d,1}^{(k)} = o(1)$ as $d \nearrow \infty$. Therefore,

$$\overline{b}_d^{(k)} = \beta_{d,0}^{(k)} + \beta_{d,1}^{(k)} = \frac{\alpha(\mathbf{r}, \mathbf{c})^k}{k!} + o(1)$$

as $d \nearrow \infty$. In order to complete the argument, recall that for each $c \in (0, \infty)$,

$$\lim_{k \longrightarrow 0} \sup_{0 \leq x \leq c} \left| \exp(x) \sum_{j=0}^{k} (-1)^j \frac{x^j}{j!} - 1 \right| = 0.$$

Under our current assumptions, we have that $\alpha(\mathbf{r}, \mathbf{c}) = O(1)$, therefore, the previous estimate for the exponential function together with (9) yields the conclusion of the result. □

**4. Complexity notions in rare-event simulation.** In these section, we shall briefly discuss basic notions of efficiency that are helpful to calculate the computational cost (in terms of the number of replications) of estimating small probabilities via simulation using an estimator of the form (4); for more on efficiency of rare-event simulation estimators [see Asmussen and Glynn (2007), Bucklew (2004) and Juneja and Shahabuddin (2006)].

Let $\beta \triangleq P(Z \in A)$ and suppose that $\beta \approx 0$. In order to be precise, we shall introduce a parameter $d$ such that $\beta_d \triangleq P(Z_d \in A) \longrightarrow 0$ as $d \nearrow \infty$ and perform our cost analysis under this asymptotic regime.

Our goal is to produce an estimator, $\widehat{\beta}_{d,k}$, with the property that for given $\varepsilon, \delta \in (0, 1)$, $|\widehat{\beta}_{d,k} - \beta_d| \leq \beta_d \varepsilon$ with probability $(1 - \delta)$. If $\widehat{\beta}_{d,k}$ has this property, we say that $\widehat{\beta}_{d,k}$ has $\varepsilon$-relative precision with $1 - \delta$ confidence. Here, we use the subindex $k$ to denote the number of i.i.d. replications required to produce $\widehat{\beta}_{d,k}$. Let $(L_{d,j} : j \geq 1)$ be i.i.d. r.v.'s such that $EL_{d,j} = \beta_d$ and consider the unbiased estimator

$$\widehat{\beta}_{d,k} = \frac{1}{k} \sum_{j=1}^{k} L_{d,j}.$$

A standard way to measure the efficiency of the estimator $\widehat{\beta}_{d,k}$ in the rare-event simulation literature relates to its variance measured in relative terms. This approach gives rise to the notion of *strong efficiency*. More precisely, if $\sigma_d^2 = \text{Var}(L_d) < \infty$, then the $L_d$'s are said to be strongly efficient if the corresponding coefficient of variation, $cv_d \triangleq \sigma_d/\beta_d$, is uniformly bounded for $d \geq 0$. In particular, in the context of importance sampling estimators discussed in Section 2.2, see (4), $L_{d,j} = L_{d,j}(Z_{d,j})$ and $\sigma_d^2 = \text{Var}_G(L_d)$. In other words, the variance must be computed according to the underlying importance sampling distribution.

One often says that $L_d$ is strongly efficient meaning that the family of $L_d$'s is strongly efficient. In order to motivate strong efficiency in terms of the computational cost (measured by the number of i.i.d. replications) required to produce an estimator that has $\varepsilon$-relative precision with $1 - \delta$ confidence, one can use Chebyshev's inequality to obtain

$$P(|\widehat{\beta}_{d,k} - \beta_d| \geq \varepsilon \beta_d) \leq \frac{\sigma_d^2}{k\varepsilon^2 \beta_d^2}.$$

Therefore, $k \geq \varepsilon^{-2}\delta^{-1}(\sigma_d/\beta_d)^2$ replications are sufficient to produce an estimator that achieves $\varepsilon$ relative precision with $1 - \delta$ confidence. Consequently, if $L_d$ is strongly efficient, the number of replications required to obtain $\varepsilon$-relative precision with $1 - \delta$ confidence is bounded as $\beta_d \longrightarrow 0$. Obviously, strong efficiency alone is not a useful concept for measuring computational complexity because nothing has been said about the computational cost attached to each replication.

When dealing with discrete structures, such as binary contingency tables, it makes sense to measure the cost per replication in terms of the amount of information (number of bits) required to encode the family of problems at hand (i.e., *the size of the problem*). In the context of binary contingency tables, statistical applications such as those described by Chen et al. (2005), require estimating the whole distribution of statistics that depend on all the entries in the table in order to perform an hypothesis test. As a consequence, it makes sense to parameterize the size of the problem, say $d$, in terms of the number of bits required to encode a binary table, which can be taken to be the number of ones (or the number zeros, but if the table is sparse, it is obviously cheaper to encode it in terms of the number of ones).

The total complexity involves multiplying the number of replications, $k$, times the cost attached to the generation of each replication which we shall denote by $\kappa(d)$ (the cost per replication is measured by the total number of operations such as additions, multiplications and comparisons in terms of the size of the problem). Therefore, in the presence of strong efficiency, by setting the number of replications $k = \Theta(\varepsilon^{-2}\delta^{-1})$, we see that $\widehat{\beta}_{d,k}$ requires $O(\kappa(d)\varepsilon^{-2}\delta^{-1})$ operations as $d \nearrow \infty$ and $\varepsilon, \delta \searrow 0$ to achieve $\varepsilon$-relative precision with $1 - \delta$ confidence.

The notions of efficiency discussed in the previous paragraph are related to standard notions found in randomized algorithms and approximate counting, such as that of *fully polynomial randomized approximation schemes* (FPRAS) [see, Mitzenmacher (2005), page 254]. In particular, an algorithm that outputs and estimator that has $\varepsilon$-relative precision with $1 - \delta$ confidence in $O(\kappa(d)\varepsilon^{-k_1} \times \log(\delta^{-1})^{k_2})$ operations, for some $k_1$, $k_2 > 0$, as $d \nearrow \infty$ and $\varepsilon, \delta \searrow 0$ is a FPRAS if $\kappa(d)$ grows polynomially in the size of the problem, say, $d$. Because of the factor $\log(\delta^{-1})^{k_2}$ that appears in the definition of a FPRAS, which is much smaller than the factor $\delta^{-1}$ that arises in the context of strong efficiently, we introduce a stronger form of efficiency that we shall call *exponential efficiency.*

DEFINITION. We say that the family of estimators $(L_d : d \geq 1)$ is *exponentially efficient* for estimating $\beta_d$ if there exists $\theta > 0$ such that

$$(12) \qquad \psi(\theta) \triangleq \sup_{d \geq 1} \log E \exp(\theta L_d / \beta_d) < \infty.$$

REMARK. In the context of importance sampling estimators introduced in Section 2.2, the expectation in (12) is taken with respect to the underlying importance sampling distribution.

The next lemma, which is a uniform version of Chernoff's bound, will be useful to relate an estimator of the form $\widehat{\beta}_{d,k}$ to a FPRAS.

LEMMA 2. *Suppose that the family of estimators $(L_d : d \geq 1)$ is exponentially efficient for estimating $\beta_d$, then for $\varepsilon > 0$ we have*

$$(13) \qquad P(|\widehat{\beta}_{d,k} - \beta_d| \geq \varepsilon \beta_d) \leq 2 \exp(-k \min(I(\varepsilon), I(-\varepsilon)),$$

*where $I(h) = \sup_\theta (\theta(1 + h) - \psi(\theta))$. Moreover, $I(\varepsilon), I(-\varepsilon) > 0$ and $I(h) \geq \rho h^2$ for some $\rho > 0$.*

PROOF. Just as in the proof of Chernoff's bound, (13) follows by an application of Chebyshev's inequality. Let $\psi_d(\theta) = \log E \exp(\theta L_d / \beta_d)$ we then obtain

$$P(\widehat{\beta}_{d,k} - \beta_d \geq \varepsilon \beta_d) \leq \exp\left(-k \sup_{\theta \geq 0}(\theta(1 + \varepsilon) - \psi_d(\theta))\right)$$

$$= \exp\left(-k \sup_\theta (\theta(1 + \varepsilon) - \psi_d(\theta))\right) \leq \exp(-kI(\varepsilon)).$$

Similarly, one obtains

$$P(\beta_d - \widehat{\beta}_{d,k} \geq \varepsilon \beta_d) \leq \exp(-kI(-\varepsilon)).$$

Inequality (13) is obtained by adding up the left- and right-hand sides of the previous displays after simple manipulations. The last part of the lemma follows from the convexity of $\psi(\cdot)$ (supremum of convex functions is convex) combined with the fact that $\psi_d(\theta) - \theta = cv_d^2\theta^2/2 + O(\theta^3)$ as $\theta \searrow 0$ uniformly over $d$ (which holds by Taylor's theorem and exponential efficiency) and the bound $\sup_{d\geq 1} cv_d^2 < \infty$ (which once again follows from exponential efficiency). $\square$

One way to verify exponential efficiency is by showing that $L_d$ is bounded above by some deterministic constant, say $c_d^*$, such that $c_d^* = O(\beta_d)$ as $d \nearrow \infty$. An immediate consequence of the previous result is that if the family $(L_d : d \geq 1)$ is exponentially efficient and $\kappa(d)$ operations are required to generate a single replication of $L_d$. Then by setting $k = \Theta(\varepsilon^{-2}\log(\delta^{-1}))$, we see that $\widehat{\beta}_{d,k}$ requires $O(\kappa(d)\varepsilon^{-2}\log(\delta^{-1}))$ operations to achieve $\varepsilon$-relative precision with $1-\delta$ confidence. If $\kappa(d)$ grows polynomially in the size of the problem $d$, then $\widehat{\beta}_{d,k}$ is the output of a FPRAS.

In addition to considerations related to the way in which the coverage parameter $\delta$ enters the complexity analysis [in the form $\delta^{-1}$ for strongly efficient estimators and $\log(\delta^{-1})$ in the context of exponential efficiency], exponential efficiency guarantees robustness properties that are desirable in practice when constructing confidence intervals via the central limit theorem [see the discussion in L'Ecuyer et al. (2008)].

**5. Complexity analysis.** This section is dedicated to the proof of the following theorem which is our main result.

THEOREM 2. *Suppose that $\max_{i\leq m} r_i = O(1)$, $\max_{j\leq n} c_j = o(d^{1/2})$ and that $[\mathbf{c}^2]_1 = O(d)$ as $d \nearrow \infty$:*

(i) *Then the estimator $L$ provided by Algorithm 1 is strongly efficient as $d \nearrow \infty$. Since according to Lemma 1, each replication of $L$ requires $O(d^3)$ operations, the computational complexity required to estimate $u(\mathbf{r}, \mathbf{c})$ with $\varepsilon$-relative precision and $(1-\delta)$ confidence is of order $O(\varepsilon^2\delta^{-1}d^3)$ as $d \nearrow \infty$ and $\varepsilon, \delta \searrow 0$.*

(ii) *Moreover, if in addition we have that $\max c_j = o(d^{1/4-\delta_0})$ for some $\delta_0 > 0$ as $d \nearrow \infty$, then the estimator $L$ provided by Algorithm 1 is exponentially efficient as $d \nearrow \infty$. Consequently, Lemma 1 implies that $O(\varepsilon^{-2}\log(\delta^{-1})d^3)$ operations are required to estimate $u(\mathbf{r}, \mathbf{c})$ with $\varepsilon$ relative precision and $(1-\delta)$ confidence as $d \nearrow \infty$ and $\varepsilon, \delta \searrow 0$.*

The following basic result (whose proof is given at the end of this section) will be very useful in the analysis of the likelihood ratio produced by our importance sampler.

LEMMA 3. *Let $\{x_j : j \geq 1\}$ be a sequence of positive integers and let us write $\{x_{i,n} : 1 \leq i \leq n\}$ to denote any nonincreasing arrangement of the set*

$\{x_i : 1 \leq i \leq n\}$ so that

$$x_{1,n} \geq x_{2,n} \geq \cdots \geq x_{n,n}.$$

Define $y_{k,n}^{(1)} = \sum_{j=k+1}^{n} x_{j,n}$ and $y_{k,n}^{(2)} = \sum_{j=k+1}^{n} x_{j,n}^2$ for $0 \leq k \leq n-1$. Then:

(i)
$$\frac{y_{k+1,n}^{(2)}}{y_{k+1,n}^{(1)}} \leq \frac{y_{k,n}^{(2)}}{y_{k,n}^{(1)}} \leq \frac{y_{0,n}^{(2)}}{y_{0,n}^{(1)}}.$$

(ii) If $y_{0,n}^{(2)}/y_{0,n}^{(1)} = O(1)$ as $n \nearrow \infty$, then there exists a constant $a > 0$ (independent of $n$ and $k$) such that

(14) $$y_{k,n}^{(2)} \leq a(n-k)$$

as $n \nearrow \infty$. Moreover, if $x_{1,n} = o(n^{\beta_0 - \delta_0})$ for $0 \leq \delta_0 < \beta_0 \leq 1/2$ then we also have that

(15) $$\frac{x_{j,n}}{y_{j-1,n}^{(1)}} \leq \frac{a^{1/2}}{1 + (n-j)^{1-\beta_0+\delta_0}}.$$

(iii) Under the assumptions of part (ii), if $\delta_0 > 0$, then

(16) $$\sup_{n \geq 1} \sum_{j=1}^{n} \frac{x_{j,n}^{1/\beta_0}}{y_{j-1,n}^2} < \infty.$$

The previous result will be applied repeatedly to the sequence of $c_k$'s which is assumed to be ordered in a nonincreasing way, namely, $c_1 \geq c_2 \geq \cdots \geq c_n$. So, for instance, assuming that $[\mathbf{c}^2] = O(d)$, then given $\boldsymbol{\rho}_0 = \mathbf{c}$ and

$$\boldsymbol{\rho}_k = (\boldsymbol{\rho}_{k,1}, \ldots, \boldsymbol{\rho}_{k,size(\boldsymbol{\rho}_k)}) = (c_{k+1}, \ldots, c_n)$$

for $j \leq n-1$, (15) implies that there exists $n_0$ such that for all $k \leq n - n_0$ we have that then $\rho_{k,1}/[\boldsymbol{\rho}_k]_1 \leq 1/2$. Similar implications are immediate from Lemma 3 and will be invoked in our future discussion.

We now proceed with the development behind Theorem 2, we first start with part (ii). By running Algorithm 1, we obtain the estimator

$$L = L_d \triangleq \prod_{k=0}^{\tau-1} \frac{w(\mathbf{S}_k, \boldsymbol{\rho}_k)}{v(\mathbf{S}_{k+1}, \boldsymbol{\rho}_{k+1})} I(\mathbf{S}_\tau = 0)$$

$$= \prod_{k=0}^{n-1} \frac{w(\mathbf{S}_k, \boldsymbol{\rho}_k)}{v(\mathbf{S}_{k+1}, \boldsymbol{\rho}_{k+1})} I(\mathbf{S}_n = 0)$$

$$= \frac{v(\mathbf{r}, \mathbf{c})}{v(\mathbf{0}, *)} \prod_{k=0}^{n-1} \frac{w(\mathbf{S}_k, \boldsymbol{\rho}_k)}{v(\mathbf{S}_k, \boldsymbol{\rho}_k)} I(\mathbf{S}_n = 0),$$

where (as indicated before) $v(\mathbf{0},*)$ is defined as 1. Recall that
$$v(\mathbf{r},\mathbf{c}) \sim u(\mathbf{r},\mathbf{c})$$
as $d \nearrow \infty$. Therefore, in order to show exponential or strong efficiency, we must study the properties of $R_d$ defined as
$$R_d((\mathbf{S}_0,\boldsymbol{\rho}_0),\ldots,(\mathbf{S}_{n-1},\boldsymbol{\rho}_{n-1})) \triangleq \prod_{k=0}^{n-1} \frac{w(\mathbf{S}_k,\boldsymbol{\rho}_k)}{v(\mathbf{S}_k,\boldsymbol{\rho}_k)} I(\mathbf{S}_n=0),$$
given $\mathbf{S}_0 = \mathbf{r}$ and $\boldsymbol{\rho}_0 = \mathbf{c}$. The analysis of $R_d$ involves studying the ratio $w(\mathbf{s}_0,\boldsymbol{\rho}_0)/v(\mathbf{s}_0,\boldsymbol{\rho}_0)$
$$\frac{w(\mathbf{s}_0,\boldsymbol{\rho}_0)}{v(\mathbf{s}_0,\boldsymbol{\rho}_0)} = \binom{m}{\rho_{0,1}}^{-1} \sum_{(\mathbf{r},\mathbf{c}) \to (\mathbf{s},\boldsymbol{\rho})} \frac{v(\mathbf{s}_1,\boldsymbol{\rho}_1)}{v(\mathbf{s}_0,\boldsymbol{\rho}_0)}.$$

Note that
$$\frac{v(\mathbf{s}_1,\boldsymbol{\rho}_1)}{v(\mathbf{s}_0,\boldsymbol{\rho}_0)} = \frac{\varphi(\mathbf{s}_1,\boldsymbol{\rho}_1)\eta(\boldsymbol{\rho}_0,m)}{\varphi(\mathbf{s}_0,\boldsymbol{\rho}_0)\eta(\boldsymbol{\rho}_1,m)}$$
and (using the notation $\rho_{i,j}$ to denote the $j$th component of the vector $\boldsymbol{\rho}_i$ for $i \in \{0,1\}$ and recalling that $\rho_{1,i} = \rho_{0,i+1}$)
$$\frac{\eta(\boldsymbol{\rho}_0,m)}{\eta(\boldsymbol{\rho}_1,m)} = \binom{m}{\rho_{0,1}}\cdots\binom{m}{\rho_{0,n}}\left(\binom{m}{\rho_{0,2}}\cdots\binom{m}{\rho_{0,n}}\right)^{-1} = \binom{m}{\rho_{0,1}}.$$

Next, observe that $\mathbf{s}_1$ is obtained by selecting a set $\Gamma = \{i_1,\ldots,i_{\rho_{0,1}}\}$ of (ordered) subindices and by picking the $i$th component of the vector $\mathbf{s}_1$, namely $s_{1,i}$, via $\mathbf{s}_{1,i} = s_{0,i} - 1$ ($i \in \Gamma$). Consequently, we have
$$\frac{\varphi(\mathbf{s}_1,\boldsymbol{\rho}_1)}{\varphi(\mathbf{s}_0,\boldsymbol{\rho}_0)} = \binom{d_0}{\rho_{0,1}}^{-1} (\Pi_{i\in\Gamma} s_{0,i}) \exp(-(\alpha(\mathbf{s}_1,\boldsymbol{\rho}_1) - \alpha(\mathbf{s}_0,\boldsymbol{\rho}_0))),$$
where $d_0 = \sum_{j=1}^m s_{0,j} = \sum_{j=1}^{n_0} \rho_{0,j}$ [with $n_0 = size(\boldsymbol{\rho}_0)$]. Therefore,
$$\frac{w(\mathbf{s}_0,\boldsymbol{\rho}_0)}{v(\mathbf{s}_0,\boldsymbol{\rho}_0)} = \binom{d_0}{\rho_{0,1}}^{-1} \sum_{(\mathbf{s}_0,\boldsymbol{\rho}_0) \to (\mathbf{s}_1,\boldsymbol{\rho}_1)} (\Pi_{i\in\Gamma} s_{0,i}) \exp(-(\alpha(\mathbf{s}_1,\boldsymbol{\rho}_1) - \alpha(\mathbf{s}_0,\boldsymbol{\rho}_0))).$$

Let us provide a more convenient expression for the previous ratio. First, we write $\mathbf{s}_{0,\Gamma} = (s_{i_1},\ldots,s_{i_{\rho_{0,1}}})$ and define $\gamma = [\boldsymbol{\rho}_1]_2/(2[\boldsymbol{\rho}_1]_1^2)$. Then (using $\mathbf{1}$ to denote the vector of ones) we have
$$\alpha(\mathbf{s}_1,\boldsymbol{\rho}_1) = \gamma[\mathbf{s}_1]_2 = \gamma([\mathbf{s}_0]_2 - [2(\mathbf{s}_{0,\Gamma} - \mathbf{1})]_1),$$
$$\alpha(\mathbf{s}_0,\boldsymbol{\rho}_0) = \left(\gamma\left(1 - \frac{\rho_{0,1}}{[\boldsymbol{\rho}_0]_1}\right)^2 + \frac{[\rho_{0,1}]_2}{2[\boldsymbol{\rho}_0]_1^2}\right)[\mathbf{s}_0]_2$$
$$= \left(\gamma - \frac{2\gamma\rho_{0,1}}{[\boldsymbol{\rho}_0]_1} + \gamma\frac{\rho_{0,1}^2}{[\boldsymbol{\rho}_0]_1^2} + \frac{[\rho_{0,1}]_2}{2[\boldsymbol{\rho}_0]_1^2}\right)[\mathbf{s}_0]_2$$

$$= \gamma[\mathbf{s}_0]_2 - \frac{2\gamma\rho_{0,1}[\mathbf{s}_0]_2}{[\boldsymbol{\rho}_0]_1} + \frac{\gamma\rho_{0,1}^2[\mathbf{s}_0]_2}{[\boldsymbol{\rho}_0]_1^2} + \frac{[\rho_{0,1}]_2[\mathbf{s}_0]_2}{2[\boldsymbol{\rho}_0]_1^2}.$$

Therefore, we have

(17)
$$\alpha(\mathbf{s}_1, \boldsymbol{\rho}_1) - \alpha(\mathbf{s}_0, \boldsymbol{\rho}_0)$$
$$= -\gamma[2(\mathbf{s}_{0,\Gamma} - \mathbf{1})]_1 + \frac{2\gamma\rho_{0,1}[\mathbf{s}_0]_2}{[\boldsymbol{\rho}_0]_1} - \frac{\gamma\rho_{0,1}^2[\mathbf{s}_0]_2}{[\boldsymbol{\rho}_0]_1^2} - \frac{[\rho_{0,1}]_2[\mathbf{s}_0]_2}{2[\boldsymbol{\rho}_0]_1^2}.$$

Define $\beta(\mathbf{s}_0, \boldsymbol{\rho}_0)$ via

$$\log \beta(\mathbf{s}_0, \boldsymbol{\rho}_0) = -\frac{2\gamma\rho_{0,1}[\mathbf{s}_0]_2}{[\boldsymbol{\rho}_0]_1} + \frac{\gamma\rho_{0,1}^2[\mathbf{s}_0]_2}{[\boldsymbol{\rho}_0]_1^2} + \frac{[\rho_{0,1}]_2[\mathbf{s}_0]_2}{2[\boldsymbol{\rho}_0]_1^2}$$

and

$$h(\mathbf{s}_{0,\Gamma}, s_{0,1}) = \Pi_{j=1}^{\rho_{0,1}}(s_{0,i_j}\exp(2\gamma(s_{0,i_j} - 1))).$$

We now are ready to provide an estimate for the ratio $w(\mathbf{s}_0, \boldsymbol{\rho}_0)/v(\mathbf{s}_0, \boldsymbol{\rho}_0)$.

LEMMA 4. *Assuming that $\max_{i\leq m} r_i = O(1)$ and that $[\mathbf{c}^2]_1 = O(d)$ as $d \nearrow \infty$ there exists a constant $\lambda \in (0, \infty)$ such that*

$$\frac{w(\mathbf{s}_0, \boldsymbol{\rho}_0)}{v(\mathbf{s}_0, \boldsymbol{\rho}_0)} \leq \exp\left(\lambda\frac{[\rho_{0,1}]_4}{[\boldsymbol{\rho}_0]_1^2}\right).$$

PROOF. Let us define $\rho_{0,1}$ i.i.d. random variables $J_1, J_2, \ldots, J_{\rho_{0,1}}$ with distribution

$$\widetilde{P}(J_1 = j) = \frac{\exp(2\gamma s_{0,j})s_{0,j}}{\widetilde{w}},$$

where $\widetilde{w} = \sum_{j=1}^{m} \exp(2\gamma s_{0,j})s_{0,j}$ and $m = size(\mathbf{s}_0)$. In addition, we shall use $\widetilde{E}(\cdot)$ to denote the expectation operator associated with $\widetilde{P}(\cdot)$ and define the event $A = \{J_i \neq J_j : i \neq j\}$ (i.e., all the $J_i$'s are different). We have

$$\binom{[\mathbf{s}_0]_1}{\rho_{0,1}}^{-1} \sum_{\Gamma \subset \{1,\ldots,m\}} h(\mathbf{s}_{0,\Gamma}, s_{0,1}) = \binom{[\mathbf{s}_0]_1}{\rho_{0,1}}^{-1} \frac{\widetilde{w}^{\rho_{0,1}}}{\rho_{0,1}!}\exp(-2\gamma\rho_{0,1})\widetilde{P}(A).$$

Let us first analyze $\widetilde{w}$. Note that under our assumptions

$$\widetilde{w} = \sum_{j=1}^{m} s_{0,j}\left(1 + 2\gamma s_{0,j} + (2\gamma)^2\frac{s_{0,j}^2}{2!} + \cdots\right)$$
$$\leq [\mathbf{s}_0]_1 \exp\left(\frac{2\gamma[\mathbf{s}_0^2]_1}{[\mathbf{s}_0]_1} + O\left(\frac{\gamma^2[\mathbf{s}_0^3]_1}{[\mathbf{s}_0]_1}\right)\right).$$

We then conclude that

$$\left(\begin{array}{c}[\mathbf{s}_0]_1\\ \rho_{0,1}\end{array}\right)^{-1}\frac{\widetilde{w}^{\rho_{0,1}}}{\rho_{0,1}!}$$

(18)
$$\leq \left(\begin{array}{c}[\mathbf{s}_0]_1\\ \rho_{0,1}\end{array}\right)^{-1}\frac{[\mathbf{s}_0]_1^{\rho_{0,1}}}{\rho_{0,1}!}\times \exp\left(\frac{2\gamma\rho_{0,1}[\mathbf{s}_0^2]_1}{[\mathbf{s}_0]_1}+O\left(\frac{\gamma^2[\mathbf{s}_0^3]_1}{[\mathbf{s}_0]_1}\right)\right)$$

$$= \prod_{k=1}^{\rho_{0,1}-1}\left(1-\frac{k}{[\mathbf{s}_0]_1}\right)^{-1}\times \exp\left(\frac{2\gamma\rho_{0,1}[\mathbf{s}_0^2]_1}{[\mathbf{s}_0]_1}+O\left(\frac{\gamma^2[\mathbf{s}_0^3]_1}{[\mathbf{s}_0]_1}\right)\right).$$

Assuming $\rho_{0,1}/[\mathbf{s}_0]_1 \leq 1/2$, we have

$$-\sum_{k=1}^{\rho_{0,1}-1}\log\left(1-\frac{k}{[\mathbf{s}_0]_1}\right) \leq \frac{[\rho_{0,1}]_2}{2[\mathbf{s}_0]_1}+\frac{\rho_{0,1}(\rho_{0,1}+1)(2\rho_{0,1}+1)}{6[\mathbf{s}_0]_1^2}$$

and, therefore, we have

$$\log\left(\left(\begin{array}{c}[\mathbf{s}_0]_1\\ \rho_{0,1}\end{array}\right)^{-1}\frac{\widetilde{w}^{\rho_{0,1}}}{\rho_{0,1}!}\right)$$

$$\leq \frac{2\gamma\rho_{0,1}[\mathbf{s}_0^2]_1}{[\mathbf{s}_0]_1}+\frac{[\rho_{0,1}]_2}{2[\mathbf{s}_0]_1}+\frac{\rho_{0,1}(\rho_{0,1}+1)(2\rho_{0,1}+1)}{6[\mathbf{s}_0]_1^2}$$

$$+O\left(\frac{\gamma^2[\mathbf{s}_0^3]_1\rho_{0,1}}{[\mathbf{s}_0]_1}\right).$$

We now estimate $\widetilde{P}(A^c)$ using the inclusion-exclusion principle and Bonferroni's inequalities. We have that

$$\widetilde{P}(A^c) \leq \binom{\rho_{0,1}}{2}\frac{1}{\widetilde{w}^2}\sum_{j=1}^{m}s_{0,j}^2\exp(4\gamma s_{0,j}),$$

this corresponds to the union bound taking all possible ways in which $\{J_{i_1} = J_{i_2}\}$ for $i_1 \neq i_2$. Next, we obtain the following lower bound corresponding to the cases in which $\{J_{i_1} = J_{i_2}, J_{i_3} = J_{i_4}\}$,

$$\widetilde{P}(A^c) \geq \binom{\rho_{0,1}}{2}\frac{1}{\widetilde{w}^2}\sum_{j=1}^{m}s_{0,j}^2\exp(4\gamma s_{0,j})$$

$$-\binom{3}{1}\binom{\rho_{0,1}}{3}\frac{1}{\widetilde{w}^3}\sum_{j=1}^{m}s_{0,j}^3\exp(12\gamma s_{0,j})$$

$$-\binom{4}{2}\binom{\rho_{0,1}}{4}\frac{1}{\widetilde{w}^4}\left(\sum_{j=1}^{m}s_{0,j}^2\exp(4\gamma s_{0,j})\right)^2.$$

We then conclude
$$\widetilde{P}(A) \leq 1 - \binom{\rho_{0,1}}{2} \frac{1}{\widetilde{w}^2} \sum_{j=1}^{m} s_{0,j}^2 \exp(4\gamma s_{0,j})$$
$$+ \binom{3}{1} \binom{\rho_{0,1}}{3} \frac{1}{\widetilde{w}^3} \sum_{j=1}^{m} s_{0,j}^3 \exp(12\gamma s_{0,j})$$
$$+ \binom{4}{2} \binom{\rho_{0,1}}{4} \frac{1}{\widetilde{w}^4} \left( \sum_{j=1}^{m} s_{0,j}^2 \exp(4\gamma s_{0,j}) \right)^2.$$

Now, we have that
$$\binom{\rho_{0,1}}{2} \frac{1}{\widetilde{w}^2} \sum_{j=1}^{m} s_{0,j}^2 \exp(4\gamma s_{0,j})$$
$$= \binom{\rho_{0,1}}{2} \frac{[\mathbf{s}_0^2]_1}{[\mathbf{s}_0]_1} \times \frac{(1 + 4\gamma[\mathbf{s}_0^3]_1/[\mathbf{s}_0^2]_1 + (4\gamma)^2[\mathbf{s}_0^4]_1/(2![\mathbf{s}_0^2]_1) + \cdots)}{(1 + 2\gamma[\mathbf{s}_0^2]_1/[\mathbf{s}_0]_1 + (2\gamma)^2[\mathbf{s}_0^3]_1/(2![\mathbf{s}_0]_1) + \cdots)}$$
$$= \binom{\rho_{0,1}}{2} \frac{[\mathbf{s}_0^2]_1}{[\mathbf{s}_0]_1} \left( 1 + 4\gamma \frac{[\mathbf{s}_0^3]_1}{[\mathbf{s}_0^2]_1} + \Theta\left( \gamma^2 \frac{[\mathbf{s}_0^4]_1}{[\mathbf{s}_0^2]_1} \right) \right)$$
$$\times \left( 1 - 2\gamma \frac{[\mathbf{s}_0^2]_1}{[\mathbf{s}_0]_1} + \Theta\left( \gamma^2 \frac{[\mathbf{s}_0^3]_1}{[\mathbf{s}_0]_1} \right) \right)$$
$$= \binom{\rho_{0,1}}{2} \frac{[\mathbf{s}_0^2]_1}{[\mathbf{s}_0]_1} \left( 1 + 4\gamma \frac{[\mathbf{s}_0^3]_1}{[\mathbf{s}_0^2]_1} - 2\gamma \frac{[\mathbf{s}_0^2]_1}{[\mathbf{s}_0]_1} + \Theta\left( \gamma^2 \frac{[\mathbf{s}_0^4]_1}{[\mathbf{s}_0]_1} \right) \right).$$

Note that
$$\frac{2[\mathbf{s}_0^3]_1}{[\mathbf{s}_0^2]_1} - \frac{[\mathbf{s}_0^2]_1}{[\mathbf{s}_0]_1}$$
$$= \frac{2[\mathbf{s}_0^3]_1[\mathbf{s}_0]_1 - [\mathbf{s}_0^2]_1^2}{[\mathbf{s}_0]_1[\mathbf{s}_0^2]_1}$$
$$= \frac{[\mathbf{s}_0^4]_1}{[\mathbf{s}_0]_1[\mathbf{s}_0^2]_1} + 2 \frac{\sum_{i<j}(s_{0,i}^3 s_{0,j} + s_{0,i} s_{0,j}^3 - s_{0,i}^2 s_{0,j}^2)}{[\mathbf{s}_0]_1[\mathbf{s}_0^2]_1}$$
$$\geq \frac{[\mathbf{s}_0^4]_1}{[\mathbf{s}_0]_1[\mathbf{s}_0^2]_1} + 2 \frac{\sum_{i<j} s_{0,i} s_{0,j} \min(s_{0,j}^2, s_{0,i}^2)}{[\mathbf{s}_0]_1[\mathbf{s}_0^2]_1}$$
$$\geq \frac{[\mathbf{s}_0^4]_1}{[\mathbf{s}_0]_1[\mathbf{s}_0^2]_1} + \frac{[\rho_{0,1}]_2}{[\mathbf{s}_0]_1[\mathbf{s}_0^2]_1}.$$

As a consequence,
$$\binom{\rho_{0,1}}{2} \frac{1}{\widetilde{w}^2} \sum_{j=1}^{m} s_{0,j}^2 \exp(4\gamma s_{0,j})$$

$$\geq \binom{\rho_{0,1}}{2} \frac{[\mathbf{s}_0^2]_1}{[\mathbf{s}_0]_1^2} \times \left(1 + \frac{2\gamma[\mathbf{s}_0^4]_1}{[\mathbf{s}_0]_1[\mathbf{s}_0^2]_1} + \frac{2\gamma[\rho_{0,1}]_2}{[\mathbf{s}_0]_1[\mathbf{s}_0^2]_1} + O\left(\gamma^2 \frac{[\mathbf{s}_0^4]_1}{[\mathbf{s}_0]_1}\right)\right)^2$$

$$\geq \binom{\rho_{0,1}}{2} \frac{[\mathbf{s}_0^2]_1}{[\mathbf{s}_0]_1^2} + \Theta\left(\frac{[\rho_{0,1}]_2^2}{[\mathbf{s}_0]_1^4} + \frac{[\rho_{0,1}]_2}{[\mathbf{s}_0]_1^3}\right).$$

Therefore,

$$\widetilde{P}(A) \leq 1 - \binom{\rho_{0,1}}{2} \frac{[\mathbf{s}_0^2]_1}{[\mathbf{s}_0]_1^2} + 3\binom{\rho_{0,1}}{3} \frac{[\mathbf{s}_0^3]_1}{[\mathbf{s}_0]_1^3}$$

$$+ \binom{4}{2}\binom{\rho_{0,1}}{4} \frac{[\mathbf{s}_0^2]_1^2}{[\mathbf{s}_0]_1^4}\left(1 + 4\gamma\frac{[\mathbf{s}_0^3]_1}{[\mathbf{s}_0^2]_1}\right)^2 + \binom{\rho_{0,1}}{2} \frac{[\mathbf{s}_0^2]_1}{[\mathbf{s}_0]_1^2}$$

$$+ O\left(\frac{[\rho_{0,1}]_2^2}{[\mathbf{s}_0]_1^4} + \frac{[\rho_{0,1}]_2}{[\mathbf{s}_0]_1^3}\right)$$

$$= 1 - \binom{\rho_{0,1}}{2} \frac{[\mathbf{s}_0^2]_1}{[\mathbf{s}_0]_1^2} + O\left(\frac{[\rho_{0,1}]_4}{[\mathbf{s}_0]_1^2} + \frac{[\rho_{0,1}]_2^2}{[\mathbf{s}_0]_1^4} + \frac{[\rho_{0,1}]_2}{[\mathbf{s}_0]_1^3}\right)$$

$$\leq \exp\left(-\binom{\rho_{0,1}}{2} \frac{[\mathbf{s}_0^2]_1}{[\mathbf{s}_0]_1^2} + O\left(\frac{[\rho_{0,1}]_4}{[\mathbf{s}_0]_1^2}\right)\right).$$

We now group all of our terms together in a convenient way in order to estimate the ratio $w(\mathbf{s}_0, \boldsymbol{\rho}_0)/v(\mathbf{s}_0, \boldsymbol{\rho}_0)$. In order to do this we define the terms $\chi_1$, $\chi_2$ and $\chi_3$ as

$$\chi_1 = -\frac{2\gamma\rho_{0,1}[\mathbf{s}_0]_2}{[\boldsymbol{\rho}_0]_1} + \frac{2\gamma\rho_{0,1}[\mathbf{s}_0^2]_1}{[\mathbf{s}_0]_1} - 2\gamma\rho_{0,1},$$

$$\chi_2 = -\binom{\rho_{0,1}}{2}\frac{[\mathbf{s}_0^2]_1}{[\mathbf{s}_0]_1^2} + \frac{[\rho_{0,1}]_2[\mathbf{s}_0]_2}{2[\boldsymbol{\rho}_0]_1^2} + \frac{[\rho_{0,1}]_2}{2[\mathbf{s}_0]_1},$$

$$\chi_3 = \frac{\rho_{0,1}(\rho_{0,1}+1)(2\rho_{0,1}+1)}{6[\mathbf{s}_0]_1^2} + \frac{\gamma\rho_{0,1}^2[\mathbf{s}_0]_2}{[\boldsymbol{\rho}_0]_1^2} + O\left(\frac{[\rho_{0,1}]_4}{[\mathbf{s}_0]_1^2} + \frac{\gamma^2[\mathbf{s}_0^3]_1\rho_{0,1}}{[\mathbf{s}_0]_1}\right).$$

Note that $\chi_3 = O([\rho_{0,1}]_4/[\mathbf{s}_0]_1^2)$, consequently if $\rho_{0,1}/[\boldsymbol{\rho}_0]_1 \leq 1/2$, then

$$\log\left(\frac{w(\mathbf{s}_0, \boldsymbol{\rho}_0)}{v(\mathbf{s}_0, \boldsymbol{\rho}_0)}\right) = \chi_1 + \chi_2 + O\left(\frac{[\rho_{0,1}]_4}{[\mathbf{s}_0]_1^2}\right).$$

Finally, we compute $\chi_1$ and $\chi_2$; first, we have that

$$-\frac{2\gamma\rho_{0,1}[\mathbf{s}_0]_2}{[\boldsymbol{\rho}_0]_1} + \frac{2\gamma\rho_{0,1}[\mathbf{s}_0^2]_1}{[\mathbf{s}_0]_1} = \frac{2\gamma\rho_{0,1}}{[\mathbf{s}_0]_1}(-[\mathbf{s}_0]_2 + [\mathbf{s}_0^2]_1)$$

$$= \frac{2\gamma\rho_{0,1}}{[\mathbf{s}_0]_1}[\mathbf{s}_0]_1 = 2\gamma\rho_{0,1}.$$

Therefore, we have that $\chi_1 = 0$. A similar computation yields $\chi_2 = 0$, indeed

$$-\binom{\rho_{0,1}}{2}\frac{[\mathbf{s}_0^2]_1}{[\mathbf{s}_0]_1^2} + \frac{[\rho_{0,1}]_2[\mathbf{s}_0]_2}{2[\boldsymbol{\rho}_0]_1^2} + \frac{[\rho_{0,1}]_2}{2[\mathbf{s}_0]_1} = \frac{[\rho_{0,1}]_2}{2[\mathbf{s}_0]_1^2}(-[\mathbf{s}_0^2]_1 + [\mathbf{s}_0]_2) + \frac{[\rho_{0,1}]_2}{2[\mathbf{s}_0]_1}$$

$$= -\frac{[\rho_{0,1}]_2}{2[\mathbf{s}_0]_1^2}[\mathbf{s}_0]_1 + \frac{[\rho_{0,1}]_2}{2[\mathbf{s}_0]_1} = 0.$$

The result of the lemma then follows. □

A consequence of the previous result is the following corollary.

COROLLARY 3. *Assuming that $\max_{i \leq m} r_i = O(1)$, $[\mathbf{c}^2]_1 = O(d)$ and $\max_{j \leq n} c_j = o(d^{1/4 - \delta_0})$ as $d \nearrow \infty$, then there exists a constant $\lambda^* \in (0, \infty)$ independent of $d$ such that*

$$R_d((\mathbf{S}_0, \boldsymbol{\rho}_0), \ldots, (\mathbf{S}_{n-1}, \boldsymbol{\rho}_{n-1})) \leq \lambda^*.$$

PROOF. Iterating the estimate obtained in Lemma 4, we obtain that

$$R_d((\mathbf{S}_0, \boldsymbol{\rho}_0), \ldots, (\mathbf{S}_{n-1}, \boldsymbol{\rho}_{n-1})) \leq \kappa_0 \exp\left(\lambda \sum_{k=0}^{n-1} \frac{[\rho_{k,1}]_4}{[\boldsymbol{\rho}_k]_1^2}\right).$$

The result then follows as a consequence of (16) in Lemma 3. □

With Corollary 3 at hand, we have all what is needed to establish exponential efficiency. However, before we put all the pieces together let us continue with the basic elements behind the strong efficiency properties indicated in part (i) of Theorem 2. We then will conclude with a summary of all our results and the complete proof of Theorem 2.

In order to establish strong efficiency we must study the function

$$g(\mathbf{s}_0, \boldsymbol{\rho}_0) \triangleq E^Q_{\mathbf{s}_0, \boldsymbol{\rho}_0}(R_d^2) = E^Q_{\mathbf{s}_0, \boldsymbol{\rho}_0}\left(\prod_{k=0}^{n-1} \frac{w^2(\mathbf{S}_k, \boldsymbol{\rho}_k)}{v^2(\mathbf{S}_k, \boldsymbol{\rho}_k)} I(\mathbf{S}_n = 0)\right).$$

In particular, we must show that $g(\mathbf{s}_0, \boldsymbol{\rho}_0)$ remains bounded as $[\mathbf{s}_0]_1 \nearrow \infty$. Our strategy is to derive a linear inequality for $g(\cdot)$ and show that one can satisfy this inequality with a convenient Lyapunov function $f(\cdot)$ that remains bounded as $d \nearrow \infty$. The next result provides sufficient conditions for the construction of an appropriate Lyapunov function. The corresponding proof is given at the end of the section.

PROPOSITION 1. *Assume $f \geq 1$ is a function that satisfies*

(19) $$f(\mathbf{s}_0, \boldsymbol{\rho}_0) \geq \frac{w^2(\mathbf{s}_0, \boldsymbol{\rho}_0)}{v^2(\mathbf{s}_0, \boldsymbol{\rho}_0)} E^Q_{\mathbf{s}_0, \boldsymbol{\rho}_0} f(\mathbf{S}_1, \boldsymbol{\rho}_1)$$

as long as $[\mathbf{s}_0]_1 \geq d_0$ [for some $d_0 \in (0, \infty)$ fixed]. Then

$$g(\mathbf{s}_0, \boldsymbol{\rho}_0) \leq \kappa_{d_0} f(\mathbf{s}_0, \boldsymbol{\rho}_0),$$

where $\kappa_{d_0} = \sup_{[\mathbf{s}_0]_1 \leq d_0} g(\mathbf{s}_0, \boldsymbol{\rho}_0) < \infty$.

A function $f(\cdot)$ satisfying the hypothesis of Proposition 1 is typically called a Lyapunov function in the context of stability of Markov chains [see Meyn and Tweedie (1993)]. Our goal is to build a bounded Lyapunov function $f(\cdot)$. Similar Lyapunov-type bounds have been studied in the rare-event simulation literature, see for instance Blanchet and Glynn (2008)] for applications in the context of rare-event estimation problems related to first passage time probabilities.

Constructing an appropriate Lyapunov function $f(\cdot)$ is typically not a simple task. Nevertheless, such construction is often guided by a solid understanding of the estimates involved in the ratio of $w(\mathbf{s}_0, \boldsymbol{\rho}_0)/v(\mathbf{s}_0, \boldsymbol{\rho}_0)$. This is precisely the strategy that we use to construction our Lyapunov function. In particular, we first put for $\theta_0 > 0$

$$f_0(\boldsymbol{\rho}) = \exp\left(\theta_0 \sum_{j=1}^{size(\boldsymbol{\rho})} \frac{\boldsymbol{\rho}_j}{(\sum_{k=j}^{n} \boldsymbol{\rho}_k)^2}\right),$$

then set

$$f_1(\mathbf{s}, \boldsymbol{\rho}) = \exp(\theta_1 \alpha(\mathbf{s}, \boldsymbol{\rho}))$$

for some $\theta_1 > 0$ and finally define

$$f(\mathbf{s}, \boldsymbol{\rho}) = f_0(\boldsymbol{\rho}) f_1(\mathbf{s}, \boldsymbol{\rho}).$$

The form of this function was obtained by inspecting carefully the analysis behind Lemma 4. We first tried a Lyapunov function such as $f_1$ and then after doing some computations, recognized the need for a term such as $f_0$ as the proof of the next lemma indicates.

LEMMA 5. *There exists $\theta_1, \theta_2 > 0$ such that $f(\cdot)$ satisfies the conditions of Proposition 1.*

PROOF. The proof proceeds along the same lines as that of Lemma 4. Given $(\mathbf{s}_0, \boldsymbol{\rho}_0)$ let us denote $(\mathbf{s}_1, \boldsymbol{\rho}_1)$ an admissible transition step [so that $(\mathbf{s}_0, \boldsymbol{\rho}_0) \to (\mathbf{s}_1, \boldsymbol{\rho}_1)$]. In particular, we have that there exists a set of subindexes $\Gamma = \{i_1, \ldots, i_{\rho_{0,1}}\}$ such that $s_{1,j} = s_{0,j} - 1(j \in \Gamma)$. We write $\gamma = [\boldsymbol{\rho}_1]_2/(2[\boldsymbol{\rho}_1]_1)$ and introduce $\rho_{0,1}$ i.i.d. random variables $J_1, \ldots, J_{\rho_{0,1}}$ with distribution

$$\widetilde{P}(J_1 = k) = \frac{\exp(2\gamma s_k) s_k}{\widetilde{w}},$$

where
$$\widetilde{w} = \sum_{i=1}^{m} \exp(2\gamma s_i) s_i.$$

We have that
$$\frac{f_1(\mathbf{s}_1, \boldsymbol{\rho}_1)}{f_1(\mathbf{s}_0, \boldsymbol{\rho}_0)} = \exp(\theta_1(\alpha(\mathbf{s}_1, \boldsymbol{\rho}_1) - \alpha(\mathbf{s}_0, \boldsymbol{\rho}_0)))$$
$$= \exp\left(\theta_1 \frac{2\gamma \rho_{0,1}[\mathbf{s}_0]_2}{[\boldsymbol{\rho}_0]_1} + 2\theta_1 \gamma \rho_{0,1} - \theta_1 \frac{\gamma \rho_{0,1}^2 [\mathbf{s}_0]_2}{[\boldsymbol{\rho}_0]_1^2}\right)$$
$$\times \exp\left(-\theta_1 \frac{[\rho_{0,1}]_2 [\mathbf{s}_0]_2}{2[\boldsymbol{\rho}_0]_1^2} - 2\theta_1 \gamma [\mathbf{s}_{0,\Gamma}]_1\right),$$

where $\mathbf{s}_{0,\Gamma} = (s_{0,j_1}, \ldots, s_{0,j_{\rho_{0,1}}})$. We need to show that there exists $d_0$, $\theta_1$ and $\theta_2$ such that for $[\mathbf{s}]_1 \geq d_0$ we have

$$\exp\left(\theta_2 \frac{\rho_{0,1}}{[\boldsymbol{\rho}_0]_1^2}\right)$$
$$\geq \frac{w(\mathbf{s}_0, \boldsymbol{\rho}_0)^2}{v(\mathbf{s}_0, \boldsymbol{\rho}_0)^2} E^Q_{\mathbf{s}, \boldsymbol{\rho}}\left(\frac{f_1(\mathbf{S}_1, \boldsymbol{\rho}_1)}{f_1(\mathbf{s}_0, \boldsymbol{\rho}_0)}\right)$$
$$= \frac{w(\mathbf{s}_0, \boldsymbol{\rho}_0)^2}{v(\mathbf{s}_0, \boldsymbol{\rho}_0)^2} \sum_{(\mathbf{s}_0, \boldsymbol{\rho}_0) \to (\mathbf{s}_1, \boldsymbol{\rho}_1)} \binom{m}{\rho_{0,1}}^{-1} \frac{f_1(\mathbf{s}_1, \boldsymbol{\rho}_1)}{f_1(\mathbf{s}_0, \boldsymbol{\rho}_0)} \frac{v(\mathbf{s}_1, \boldsymbol{\rho}_1)}{w(\mathbf{s}_0, \boldsymbol{\rho}_0)}$$
$$= \frac{w(\mathbf{s}_0, \boldsymbol{\rho}_0)}{v(\mathbf{s}_0, \boldsymbol{\rho}_0)} \sum_{(\mathbf{s}_0, \boldsymbol{\rho}_0) \to (\mathbf{s}_1, \boldsymbol{\rho}_1)} \binom{m}{\rho_{0,1}}^{-1} \frac{f_1(\mathbf{s}_1, \boldsymbol{\rho}_1)}{f_1(\mathbf{s}_0, \boldsymbol{\rho}_0)} \frac{v(\mathbf{s}_1, \boldsymbol{\rho}_1)}{v(\mathbf{s}_0, \boldsymbol{\rho}_0)}.$$

As in Lemma 4, we have that
$$\sum_{(\mathbf{s}_0, \boldsymbol{\rho}_0) \to (\mathbf{s}_1, \boldsymbol{\rho}_1)} \binom{m}{\rho_{0,1}}^{-1} \frac{f_1(\mathbf{s}_1, \boldsymbol{\rho}_1)}{f_1(\mathbf{s}_0, \boldsymbol{\rho}_0)} \frac{v(\mathbf{s}_1, \boldsymbol{\rho}_1)}{v(\mathbf{s}_0, \boldsymbol{\rho}_0)}$$
$$= \binom{[\mathbf{s}_0]_1}{\rho_{0,1}}^{-1} \widetilde{w}^{\rho_{0,1}} \exp\left((\theta_1 - 1)\frac{2\gamma \rho_{0,1}[\mathbf{s}_0]_2}{[\boldsymbol{\rho}_0]_1} + 2(\theta_1 - 1)\gamma \rho_{0,1}\right)$$
$$\times \exp\left(-(\theta_1 - 1)\frac{\gamma \rho_{0,1}^2 [\mathbf{s}_0]_2}{[\boldsymbol{\rho}_0]_1^2} - (\theta_1 - 1)\frac{[\rho_{0,1}]_2 [\mathbf{s}_0]_2}{2[\boldsymbol{\rho}_0]_1^2}\right)$$
$$\times \widetilde{E}\left(\exp\left(\sum_{i=1}^{\rho_{0,1}} -2\theta \gamma s_{0, J_i}\right); A\right),$$

where $A$ is the event that consists that all the $J_i$'s are distinct. During the proof of Lemma 4, we obtained that if $\rho_{0,1}/[\mathbf{s}_0]_1 \leq 1/2$ then
$$\log\left(\binom{[\mathbf{s}_0]_1}{\rho_{0,1}}^{-1} \frac{\widetilde{w}^{\rho_{0,1}}}{\rho_{0,1}!}\right)$$

$$\leq \frac{2\gamma\rho_{0,1}[\mathbf{s}_0^2]_1}{[\mathbf{s}_0]_1} + \frac{[\rho_{0,1}]_2}{[\mathbf{s}_0]_1} + O\left(\frac{[\rho_{0,1}]_4}{[\boldsymbol{\rho}_0]_1^2}\right).$$

Now, evidently we have that

$$\widetilde{E}\left(\exp\left(\sum_{i=1}^{\rho_{0,1}} -2\theta_1\gamma s_{0,J_i}\right); A\right) \leq (\widetilde{E}(\exp(-2\theta_1\gamma s_{0,J_i})))^{\rho_{0,1}}$$

$$\leq \exp(-2\theta_1\gamma\rho_{0,1}\widetilde{E}s_{0,J_i} + O(\theta_1\rho_{0,1}\gamma^2)).$$

Therefore, combining all these estimates together with Lemma 4 we have that there exists a constant $\lambda > 0$ such that

$$\frac{w(\mathbf{s}_0, \boldsymbol{\rho}_0)^2}{v(\mathbf{s}_0, \boldsymbol{\rho}_0)^2} E^Q_{\mathbf{s},\boldsymbol{\rho}}\left(\frac{f_1(\mathbf{S}_1, \boldsymbol{\rho}_1)}{f_1(\mathbf{s}_0, \boldsymbol{\rho}_0)}\right)$$

$$\leq \exp\left(\lambda\frac{\rho_{0,1}^4}{[\boldsymbol{\rho}_0]_1^2} + \frac{2\gamma\rho_{0,1}[\mathbf{s}_0^2]_1}{[\mathbf{s}_0]_1} + \frac{[\rho_{0,1}]_2}{[\mathbf{s}_0]_1}\right)$$

(20a)
$$\times \exp\left((\theta_1 - 1)\frac{2\gamma\rho_{0,1}[\mathbf{s}_0]_2}{[\boldsymbol{\rho}_0]_1} + 2(\theta_1 - 1)\gamma\rho_{0,1}\right)$$

$$\times \exp\left(-(\theta_1 - 1)\frac{\gamma\rho_{0,1}^2[\mathbf{s}_0]_2}{[\boldsymbol{\rho}_0]_1^2} - (\theta_1 - 1)\frac{[\rho_{0,1}]_2[\mathbf{s}_0]_2}{2[\boldsymbol{\rho}_0]_1^2}\right)$$

$$\times \exp\left(-\frac{2\theta_1\gamma\rho_{0,1}[\mathbf{s}_0^2]}{[\mathbf{s}_0]_1} + O(\theta_1\rho_{0,1}\gamma^2)\right).$$

Note that in the last line of the previous display we have used the fact that

$$\widetilde{E}s_{J_1}^2 = \frac{[\mathbf{s}_0^2]}{[\mathbf{s}_0]_1} + O(\gamma).$$

Now we note (just as we did in Lemma 4) that

$$-\frac{2\gamma\rho_{0,1}[\mathbf{s}_0]_2}{[\boldsymbol{\rho}_0]_1} + \frac{2\gamma\rho_{0,1}[\mathbf{s}_0^2]_1}{[\mathbf{s}_0]_1} - 2\gamma\rho_{0,1} = 0,$$

which implies that the logarithm of the right-hand side of (20a) equals

$$\lambda\frac{[\rho_{0,1}]_4}{[\boldsymbol{\rho}_0]_1^2} + \frac{[\rho_{0,1}]_2}{[\mathbf{s}_0]_1} - (\theta_1 - 1)\frac{\gamma\rho_{0,1}^2[\mathbf{s}_0]_2}{[\boldsymbol{\rho}_0]_1^2}$$

$$- (\theta_1 - 1)\frac{[\rho_{0,1}]_2[\mathbf{s}_0]_2}{2[\boldsymbol{\rho}_0]_1^2} + O(\theta_1\rho_{0,1}\gamma^2).$$

It is immediate from the previous expression that one can select first $\theta_1 > 0$ and then $\theta_2$ depending on $\theta_1$ so that the previous quantity is less or equal to $\theta_2\rho_{0,1}/[\boldsymbol{\rho}_0]_1^2$ as long as $[\mathbf{s}_0]_1 \geq d_0$ so that $\rho_{0,1}/[\mathbf{s}_0]_1 \leq 1/2$. The conclusion of the lemma then follows. $\square$

It is time to summarize all the previous estimates and to complete the proof of Theorem 2.

PROOF OF THEOREM 2. We first establish part (ii). By virtue of Theorem 1 and Corollary 3, we have that

$$\frac{L_d}{u(\mathbf{r},\mathbf{c})} = \frac{v(\mathbf{r},\mathbf{c})}{u(\mathbf{r},\mathbf{c})} R_d \leq \frac{v(\mathbf{r},\mathbf{c})}{u(\mathbf{r},\mathbf{c})} \lambda^* = O(1)$$

as $d \nearrow \infty$. Therefore, because of our observations in Section 4, $L_d$ is exponentially efficient and part (ii) follows. Part (i) is established similarly thanks to Lemma 5. Note that

$$\frac{E^Q_{\mathbf{r},\mathbf{c}} L_d^2}{u(\mathbf{r},\mathbf{c})^2} = \frac{v(\mathbf{r},\mathbf{c})^2}{u(\mathbf{r},\mathbf{c})^2} g(\mathbf{r},\mathbf{c}) \leq \frac{v(\mathbf{r},\mathbf{c})^2}{u(\mathbf{r},\mathbf{c})^2} f(\mathbf{r},\mathbf{c}).$$

Theorem 1 guarantees that $v(\mathbf{r},\mathbf{c})^2/u(\mathbf{r},\mathbf{c})^2 \longrightarrow 1$ as $d \nearrow \infty$. On the other hand, Lemma 3 implies that $f(\mathbf{r},\mathbf{c}) = O(1)$ as $d \nearrow \infty$. This concludes the proof of Theorem 2. □

Finally, before providing the proof of the pending results, it is worth discussing the practical implications of the previous bounds. The previous results imply a bound on the coefficient of variation that involves an exponential function to a power that depends on the maximum degree of the row sums. In practical situations, this bound can quickly become large, so the bounds given here, although computable, may be far too pessimistic in practical applications. Improving these bounds is particularly interesting given that empirically according to Chen et al. (2005), the estimated coefficient of variation of the estimator given by Algorithm 1 is consistently small (they report values that are even less than 1). The key issue involves controlling the behavior of the row sums during the course of the algorithm under $Q(\cdot)$. The techniques here can be adapted to deal with situations when the row sums may grow and this will be illustrated elsewhere in the future.

PROOF OF LEMMA 3. We have that

$$\frac{y^{(2)}_{k+1,n}}{y^{(1)}_{k+1,n}} - \frac{y^{(2)}_{k,n}}{y^{(1)}_{k,n}} = \frac{y^{(2)}_{k+1,n} y^{(1)}_{k,n} - y^{(2)}_{k,n} y^{(1)}_{k+1,n}}{y^{(1)}_{k+1,n} y^{(1)}_{k,n}}.$$

Now

$$y^{(2)}_{k+1,n} y^{(1)}_{k,n} - y^{(2)}_{k,n} y^{(1)}_{k+1,n} = (y^{(2)}_{k,n} - x^2_{k+1,n}) y^{(1)}_{k,n} - y^{(2)}_{k,n}(y^{(1)}_{k,n} - x_{k+1,n})$$

$$= x_{k+1,n}(y^{(2)}_{k,n} - x_{k+1,n} y^{(1)}_{k,n}).$$

The result then follows from the fact that

$$y_{k,n}^{(2)} = \sum_{j=k+1}^{n} x_{j,n}^2 \leq x_{k+1,n} \sum_{j=k+1}^{n} x_{j,n} = x_{k+1,n} y_{k,n}^{(1)}.$$

For part (ii) we note that, by assumption there exists $a > 0$ such that $y_{0,n}^{(2)} \leq a^{1/2} y_{0,n}^{(1)}$. Using Cauchy–Schwarz inequality and part (i) it follows that

$$y_{k,n}^{(2)} \leq a^{1/2} y_{k,n}^{(1)} \leq a^{1/2}((n-k) y_{k,n}^{(2)})^{1/2},$$

which implies $y_{k,n}^{(2)} \leq a(n-k)$. Finally, combining part (i) and the assumption that $y_{0,n}^{(2)}/y_{0,n}^{(1)} = O(1)$, we can write

$$\frac{x_{j,n}}{y_{j-1,n}^{(1)}} \leq a^{1/2} \frac{x_{j,n}}{y_{j-1,n}^{(2)}} = \frac{a^{1/2}}{x_{j,n} + x_{j+1,n}/x_{j,n} + \cdots + x_{n,n}/x_{j,n}}.$$

Now, it follows that

$$x_{j,n} + x_{j+1,n}/x_{j,n} + \cdots + x_{n,n}/x_{j,n} \geq 1 + (n-j)/x_{1,n-j}.$$

We conclude that

$$\frac{x_{j,n}}{y_{j-1,n}^{(1)}} \leq \frac{a^{1/2}}{1 + (n-j)^{1-\beta_0+\delta_0}},$$

which yields (15).

For part (iii), we use (14) and (15). In particular, we have that

$$x_{j,n} \leq \frac{a^{1/2} y_{j-1,n}^{(1)}}{1 + (n-j)^{1-\beta_0+\delta_0}} \leq a(n-j)^{\beta_0-\delta_0}$$

and, therefore,

$$\sum_{j=1}^{n} \frac{x_{j,n}^{1/\beta_0}}{(y_{j-1,n}^{(1)})^2} = O\left(\sum_{j=1}^{n-1} \frac{1}{(n-j)^{1+\delta_0/\beta_0}}\right),$$

which yields (16). □

PROOF OF PROPOSITION 1. Define $\tau_{d_0} = \inf\{k \geq 0 : \mathbf{S}_k < d_0\}$ and let $\mathcal{F}_k = \sigma(\mathbf{S}_j : 0 \leq j \leq k)$ be the $\sigma$-field generated by the process $\mathbf{S}$ up to time $k$. As in Section 3, we let $\tau$ be the first time $k \leq n$ for which the number of strictly positive components of the vector $\mathbf{S}_k$ is less than $c_{k+1}$ (and set

$c_{n+1} = 1$). Note that (using the notation $a \wedge b$ for the minimum between $a$ and $b$)

$$g(\mathbf{s}_0, \boldsymbol{\rho}_0) = E^Q_{\mathbf{s}_0, \boldsymbol{\rho}_0} \left( \prod_{k=0}^{n-1} \frac{w^2(\mathbf{S}_k, \boldsymbol{\rho}_k)}{v^2(\mathbf{S}_k, \boldsymbol{\rho}_k)} I(\mathbf{S}_n = 0) \right)$$

$$= E^Q_{\mathbf{s}_0, \boldsymbol{\rho}_0} \left( \prod_{k=0}^{\tau_{d_0} \wedge \tau - 1} \frac{w^2(\mathbf{S}_k, \boldsymbol{\rho}_k)}{v^2(\mathbf{S}_k, \boldsymbol{\rho}_k)} g(\mathbf{S}_{\tau_d \wedge \tau}, \boldsymbol{\rho}_{\tau_d \wedge \tau}) \right).$$

Clearly, the dynamics of $Q(\cdot)$ imply

$$g(\mathbf{S}_{\tau_{d_0} \wedge \tau}, \boldsymbol{\rho}_{\tau_{d_0} \wedge \tau}) 1(\tau_{d_0} > \tau) = 0,$$

therefore,

$$g(\mathbf{s}_0, \boldsymbol{\rho}_0) = E^Q_{\mathbf{s}_0, \boldsymbol{\rho}_0} \left( \prod_{k=0}^{\tau_{d_0} - 1} \frac{w^2(\mathbf{S}_k, \boldsymbol{\rho}_k)}{v^2(\mathbf{S}_k, \boldsymbol{\rho}_k)} g(\mathbf{S}_{\tau_d}, \boldsymbol{\rho}_{\tau_d}); \tau_d < \tau \right)$$

(21)

$$\leq \left( \sup_{[\mathbf{s}_0]_1 \leq d_0} g(\mathbf{s}_0, \boldsymbol{\rho}_0) \right) E^Q_{\mathbf{s}_0, \boldsymbol{\rho}_0} \left( \prod_{k=0}^{\tau_{d_0} - 1} \frac{w^2(\mathbf{S}_k, \boldsymbol{\rho}_k)}{v^2(\mathbf{S}_k, \boldsymbol{\rho}_k)} \right).$$

Now, define the stochastic process $(Z_k : k \geq 0)$ via

$$Z_k = f(\mathbf{S}_k, \boldsymbol{\rho}_k) \prod_{j=0}^{k-1} \frac{w^2(\mathbf{S}_j, \boldsymbol{\rho}_j)}{v^2(\mathbf{S}_j, \boldsymbol{\rho}_j)}$$

and consider the stopped process $M_k = Z_{k \wedge \tau_{d_0}}$. Note that $(M_k : k \geq 0)$ is a nonnegative supermartingale, that is,

$$E^Q(M_{k+1} | \mathcal{F}_k)$$
$$= E(M_{k+1}; \tau_{d_0} > k | \mathcal{F}_k) + E(M_{k+1}; \tau_{d_0} \leq k | \mathcal{F}_k)$$
$$= 1(\tau_{d_0} > k) \prod_{j=0}^{k} \frac{w^2(\mathbf{S}_j, \boldsymbol{\rho}_j)}{v^2(\mathbf{S}_j, \boldsymbol{\rho}_j)} E(f(\mathbf{S}_{k+1}, \boldsymbol{\rho}_{k+1}) | \mathbf{S}_k)$$
$$+ 1(\tau_{d_0} \leq k) \prod_{j=0}^{\tau_{d_0} - 1} \frac{w^2(\mathbf{S}_j, \boldsymbol{\rho}_j)}{v^2(\mathbf{S}_j, \boldsymbol{\rho}_j)} f(\mathbf{S}_{\tau_{d_0}}, \boldsymbol{\rho}_{k+1})$$
$$\leq 1(\tau_{d_0} > k) \prod_{j=0}^{k-1} \frac{w^2(\mathbf{S}_j, \boldsymbol{\rho}_j)}{v^2(\mathbf{S}_j, \boldsymbol{\rho}_j)} f(\mathbf{S}_k, \boldsymbol{\rho}_k)$$
$$+ 1(\tau_{d_0} \leq k) \prod_{j=0}^{\tau_{d_0} - 1} \frac{w^2(\mathbf{S}_j, \boldsymbol{\rho}_j)}{v^2(\mathbf{S}_j, \boldsymbol{\rho}_j)} f(\mathbf{S}_{\tau_{d_0}}, \boldsymbol{\rho}_{k+1})$$
$$= M_k.$$

34    J. H. BLANCHET

Therefore,

$$f(\mathbf{s}_0, \boldsymbol{\rho}_0) \geq E^Q_{\mathbf{s}_0, \boldsymbol{\rho}_0}\left( f(\mathbf{S}_{\tau_{d_0}}, \rho_{\tau_{d_0}}) \prod_{j=0}^{\tau_{d_0}-1} \frac{w^2(\mathbf{S}_j, \boldsymbol{\rho}_j)}{v^2(\mathbf{S}_j, \boldsymbol{\rho}_j)} \right)$$

$$\geq E^Q_{\mathbf{s}_0, \boldsymbol{\rho}_0}\left( \prod_{j=0}^{\tau_{d_0}-1} \frac{w^2(\mathbf{S}_j, \boldsymbol{\rho}_j)}{v^2(\mathbf{S}_j, \boldsymbol{\rho}_j)} \right).$$

This estimate, together with (21), implies the conclusion of the proposition. □

**Acknowledgments.** The author is grateful to Persi Diaconis and Jun S. Liu for helpful discussions on contingency tables, to Dan Rudoy for valuable conversations on approximate counting, to Alistair Sinclair and Alexandre Stauffer for their very insightful comments on an earlier version of this paper, and to the referee for his careful reading and suggestions that helped improve the presentation.## REFERENCES

ASMUSSEN, S. and GLYNN, P. W. (2007). *Stochastic Simulation: Algorithms and Analysis. Stochastic Modelling and Applied Probability* **57**. Springer, New York. MR2331321

BAYATI, M., KIM, J. and SABERI, A. (2007). *A Sequential Algorithm for Generating Random Graphs. Lecture Notes in Computer Science* **4627**. 326–340. Springer, Berlin.

BÉKÉSSY, A., BÉKÉSSY, P. and KOMLÓS, J. (1972). Asymptotic enumeration of regular matrices. *Studia Sci. Math. Hungar.* **7** 343–353. MR0335342

BEZÁKOVÁ, I., BHATNAGAR, N. and VIGODA, E. (2006). Sampling binary contingency tables with a greedy start. In *Proceedings of the Seventeenth Annual ACM-SIAM Symposium on Discrete Algorithms* 414–423. ACM, New York. MR2368838

BEZÁKOVÁ, I., SINCLAIR, A., ŠTEFANKOVIČ, D. and VIGODA, E. (2007). Negative examples for sequential importance sampling of binary contingency tables. In *Algorithms—ESA 2006. Lecture Notes in Computer Science* **4168** 136–147. Springer, Berlin. MR2347138

BLANCHET, J. and GLYNN, P. (2008). Efficient rare-event simulation for the maximum of heavy-tailed random walks. *Ann. Appl. Probab.* **18** 1351–1378. MR2434174

BLANCHET, J. and LIU, J. C. (2008). State-dependent importance sampling for regularly varying random walks. *Adv. in Appl. Probab.* **40** 1104–1128.

BLITZSTEIN, J. and DIACONIS, P. (2008). A sequential importance sampling algorithm for generating random graphs with prescribed degrees. Preprint.

BOTEV, Z. I. and KROESE, D. P. (2008). Non-asymptotic bandwidth selection for density estimation of discrete data. *Methodol. Comput. Appl. Probab.* **10** 435–451. MR2415127

BUCKLEW, J. A. (2004). *Introduction to Rare Event Simulation*. Springer, New York. MR2045385

CHEN, S. X. and LIU, J. S. (1997). Statistical applications of the Poisson-binomial and conditional Bernoulli distributions. *Statist. Sinica* **7** 875–892. MR1488647

CHEN, X.-H., DEMPSTER, A. P. and LIU, J. S. (1994). Weighted finite population sampling to maximize entropy. *Biometrika* **81** 457–469. MR1311090

IEOR Department
Columbia University
340 S. W. Mudd Building
500 West 120th Street
New York, New York 10027-6699
USA
E-mail: blanchet@stat.harvard.edu